\documentclass [a4paper, 12pt, reqno]{amsart}
\usepackage [latin1]{inputenc}
\usepackage {a4}
\usepackage{amscd}
\usepackage{epsfig}
\usepackage{amssymb}
\usepackage{amsmath}
\usepackage{amsthm}
\usepackage[T1]{fontenc}
\usepackage{ae,aecompl}
\usepackage{color}

\usepackage{geometry}
\newcommand{\R} {\ensuremath{\mathbb{R}}}
\newcommand{\N} {\ensuremath{\mathbb{N}}}
\newcommand{\C} {\ensuremath{\mathbb{C}}}

\newcommand{\OO}{\mathcal{O}}
\renewcommand{\o}[1]{\overline{#1}}

\newcommand{\dq}{\overline{\partial}}
\newcommand{\wt}[1]{\widetilde{#1}}
\DeclareMathOperator{\Reg}{Reg}
\DeclareMathOperator{\Sing}{Sing}

\DeclareMathOperator{\Dom}{Dom}

\newtheorem {satz} {Satz} [section]
\newtheorem {lem} [satz] {Lemma}

\newtheorem {prop} [satz] {Proposition}

\newtheorem {thm} [satz] {Theorem}

\DeclareMathOperator{\supp}{supp}

\renewcommand{\Im}{\mbox{Im }}
\newcommand{\Ra}{\mathcal{R}}

\renewcommand{\theta}{\vartheta}

\title[$\dq$ and $\dq$-Neumann operator at isolated singularities] 
{$L^2$-Properties of the $\dq$ and the $\dq$-Neumann operator on spaces with isolated singularities}

\author{N. {\O}vrelid}
\address{Department of Mathematics, University of Oslo, P.B 1053 Blindern, Oslo, N-0316 Norway.}
\email{nilsov@math.uio.no}

\author{J. Ruppenthal}
\address{Department of Mathematics, University of Wuppertal, Gau{\ss}str. 20, 42119 Wuppertal, Germany.}
\email{ruppenthal@uni-wuppertal.de}

\date{\today}



\begin{document}

\begin{abstract}
Let $X$ be a Hermitian complex space of pure dimension
with only isolated singularities and $\pi: M\rightarrow X$ a resolution of singularities.
Let $\Omega\subset\subset X$ be a domain with no singularities in the boundary,
$\Omega^*=\Omega\setminus\Sing X$ and $\Omega'=\pi^{-1}(\Omega)$.
We relate $L^2$-properties of the $\dq$ and the $\dq$-Neumann operator on $\Omega^*$
to properties of the corresponding
operators on $\Omega'$ (where the situation is classically well understood). 
Outside some middle degrees, there are compact solution operators
for the $\dq$-equation on $\Omega^*$ exactly if there are such operators on the resolution $\Omega'$,
and the $\dq$-Neumann operator is compact on $\Omega^*$ exactly if it is compact on $\Omega'$.
\end{abstract}

\maketitle

\section{Introduction}

The Cauchy-Riemann operator $\dq$ and the related $\dq$-Neumann operator play a central role in complex analysis.
Especially the $L^2$-theory for these operators is of particular importance and has become indispensable for the subject
after the fundamental work of
H\"ormander on $L^2$-estimates and existence theorems for the $\dq$-operator (see \cite{Hoe1} and \cite{Hoe2})
and the related work of Andreotti and Vesentini (see \cite{AnVe}).
By no means less important is Kohn's solution of the $\dq$-Neumann problem (see \cite{Ko1}, \cite{Ko2} and also \cite{KoNi}),
which implies existence and regularity results for the $\dq$-complex, as well (see Chapter III.1 in \cite{FK}).
Important applications of the $L^2$-theory are for instance the Ohsawa-Takegoshi extension theorem \cite{OT},
Siu's analyticity of the level sets of Lelong numbers \cite{Siu0}
or the invariance of plurigenera \cite{Siu}.

Whereas the theory of the $\dq$-operator and the $\dq$-Neumann operator is very well developed on
complex manifolds, not too much is known about the situation on singular complex spaces which
appear naturally as the zero sets of holomorphic functions.
The further development of this theory is an important endeavor
since analytic methods have led to fundamental advances in geometry on complex manifolds (see Siu's results mentioned above),
but these analytic tools are still missing on singular spaces.

The topic has attracted some attention recently and considerable progress has been made.
Let us mention e.g. the development of some Koppelman formulas by Andersson and Samuelsson (\cite{AS1}, \cite{AS2}).
Concerning the $L^2$-theory for the $\dq$-operator, {\O}vrelid and Vassiliadou
obtained essential results for the case of isolated singularities.
Following a path prepared by Pardon and Stern (see \cite{P}, \cite{PS1}, \cite{PS2}),
Forn{\ae}ss, Diederich, Vassiliadou and {\O}vrelid
(see \cite{Fo}, \cite{DFV}, \cite{FOV2}, \cite{OV1}, \cite{OV2}) and by
Ruppenthal and Zeron (see \cite{Rp1}, \cite{Rp7}, \cite{Rp8}, \cite{RZ1}, \cite{RZ2}),
they were finally able to describe the $L^2$-cohomology for the $\dq$-operator
at isolated singularities completely in terms of a resolution of singularities (see \cite{OV3}).
For another, different approach to these results we refer also to \cite{Rp10}.

This is an important progress concerning the understanding of the obstructions
to solving the $\dq$-equation at isolated singularities.
It remains to study the regularity of the equation:
On domains in complex manifolds, the close connection between the regularity of the $\dq$-equation on one hand
and the geometry of the domain (and its boundary) on the other hand is one of the central topics
of complex analysis. It is an interesting task to establish such connections also
between the regularity of the $\dq$-equation at singularities and the
geometry of the singularities. In the present paper, we study the existence
of compact solution operators for the $\dq$-equation at isolated singularities and compactness
of the $\dq$-Neumann operator in the presence of isolated singularities.
This complements the discussion of the topic in \cite{Rp9}.

Compactness can be seen as a boundary case of subelliptic regularity
(when the gain in the subelliptic estimate tends to zero), and is an important property
in the study of weakly pseudoconvex domains (see \cite{S} for a comprehensive discussion of the topic).
Moreover, compactness of the $\dq$-Neumann operator yields that the corresponding space of $L^2$-forms
has an orthonormal basis consisting of eigenforms of the $\dq$-Laplacian $\Box=\dq\dq^* +\dq^*\dq$.
The eigenvalues of $\Box$ are non-negative, have no finite limit point
and appear with finite multiplicity. It might be an interesting question to study
whether there is a nice connection between the eigenvalues and the structure of
the singularities.

\medskip

Our main results are as follows.
Let $X$ be a Hermitian complex space\footnote{%A reduced complex space with a Hermitian metric on the regular part
%which is induced by local embeddings into complex number space, hence extends continuously into the singular set.More precisely, 
A Hermitian complex space $(X,g)$ is a reduced complex space $X$ with a metric $g$ on the regular part
such that the following holds: If $x\in X$ is an arbitrary point there exists a neighborhood $U=U(x)$ and a
biholomorphic embedding of $U$ into a domain $G$ in $\C^N$ and an ordinary smooth Hermitian metric in $G$
whose restriction to $U$ is $g|_U$.}
of pure dimension $n$ with only isolated singularities. Let $\pi: M \rightarrow X$
be a resolution of singularities which exists due to Hironaka (see Section \ref{sec:resolution}), and let
$\sigma$ be any (positive definite) Hermitian metric on $M$.
We denote by $L^{p,q}$ the spaces of $L^2$-forms on $\Reg X$,
and by $L^{p,q}_\sigma$ the spaces of $L^2$-forms on $M$ with respect to $\sigma$.
Let $\Omega\subset\subset X$ be a relatively compact open subset of $X$ such that the boundary of $\Omega$
does not intersect the singular set of $X$, $b\Omega\cap \Sing X=\emptyset$.
Let $\Omega^*:=\Omega\setminus\Sing X$ and $\Omega':=\pi^{-1}(\Omega)$.

Thus, the resolution of singularities has the following nice effect:
If the original domain $\Omega$ has a "good" boundary $b\Omega$,
then $\Omega'$ is a domain in a complex manifold with the same "good" boundary.
One might consider for example a domain $\Omega$ with a strongly pseudoconvex boundary,
or assume that $X$ is a compact space and $\Omega=X$ (no boundary at all).
In both cases we know that the $\dq$-equation has compact solution operators on $\Omega'$
(modulo the obstructions to solving the equation), and that the $\dq$-Neumann operator
exists and is compact. 
It is thus interesting to relate properties of the $\dq$-operator on $\Omega^*$ (which have to be studied)
to properties of the $\dq$-operator on $\Omega'$ (which are well understood):

\begin{thm}\label{thm:1.1}
Let $q\geq 1$ and either $p+q\neq n$ or $(p,q)=(0,n)$. Under the assumptions above,
the $\dq$-operator in the sense of distributions
\begin{eqnarray*}
\dq: L^{p,q-1}(\Omega^*) \rightarrow L^{p,q}(\Omega^*)
\end{eqnarray*}
has closed range (of finite codimension) in $\ker\dq\subset L^{p,q}(\Omega^*)$
exactly if the $\dq$-operator in the sense of distributions
\begin{eqnarray*}
\dq_M: L^{p,q-1}_\sigma(\Omega') \rightarrow L^{p,q}_\sigma(\Omega')
\end{eqnarray*}
has closed range (of finite codimension) in $\ker\dq_M \subset L^{p,q}_\sigma(\Omega')$.

If this is the case,
then there exists a compact $\dq$-solution operator
\begin{eqnarray*}
{\bf S}: \Im \dq \subset L^{p,q}(\Omega^*) \rightarrow L^{p,q-1}(\Omega^*)
\end{eqnarray*}
exactly if there exists a compact $\dq$-solution operator
\begin{eqnarray*}
{\bf S}_M: \Im\dq_M \subset L^{p,q}_\sigma(\Omega') \rightarrow L^{p,q-1}_\sigma(\Omega').
\end{eqnarray*}
\end{thm}

The phrase 'of finite codimension' is optional and may or may not be included in the statement of the theorem,
just as desired. Note that there are interesting cases where $\Im \dq$ is closed without being of finite codimension in $\ker\dq$.
Consider the following example:
Let $X$ be a variety of dimension $n\geq 3$ with isolated singularities in $\C^N$
and $\Omega$ the intersection of $X$ with a spherical shell $S=\{z: r<|z-z_0|<R\}$.
Assume that $bS$ intersects $X$ transversally in $\Reg X$. 
When $\Omega'$ is a desingularization of $\Omega$, it follows from the arguments in \cite{Sh} and Theorem 3.1
in \cite{Hoe3} that the image of $\dq: L^{0,n-2}(\Omega')\rightarrow L^{0,n-1}(\Omega')$
is closed of infinite codimension in $\ker \dq \subset L^{0,n-1}(\Omega')$.
Moreover, the minimal solution operator is compact.

Our main tools in the proof of Theorem \ref{thm:1.1} are the existence of $\dq$-solution operators with some gain of regularity
at isolated singularities from \cite{FOV2} (see Theorem \ref{thm:fov2}) and a characterization of precompactness
in the space of $L^2$-forms on arbitrary Hermitian manifolds (Theorem \ref{thm:precompact2}),
which shows that these $\dq$-solution operators are compact (Theorem \ref{thm:sop}).
Other ingredients are Hironaka's resolution of singularities and Kohn's subelliptic estimates.

\smallskip

Using Theorem \ref{thm:1.1} and a representation of the $\dq$-Neumann operator in terms
of canonical solution operators for the $\dq$-equation (Theorem \ref{thm:N100}), we deduce:

\begin{thm}\label{thm:1.2}
Let $X$ be a Hermitian complex space of pure dimension $n$ with only isolated singularities,
and $\Omega\subset\subset X$ with no singularities in the boundary $b\Omega$.
Let $\pi: M\rightarrow X$ be a resolution of singularities
and $\sigma$ any (positive definite) Hermitian metric on $M$.
Let $q\geq 1$ and assume that the $\dq$-operators in the sense of distributions
with respect to the metric $\sigma$ on $\Omega'=\pi^{-1}(\Omega)$,
\begin{eqnarray*}
&\dq^M_{p,q}:& L^{p,q-1}_\sigma(\Omega')\rightarrow L^{p,q}_\sigma(\Omega'),\\
&\dq^M_{p,q+1}:& L^{p,q}_\sigma(\Omega') \rightarrow L^{p,q+1}_\sigma(\Omega'),
\end{eqnarray*}
both have closed range in the corresponding kernels of $\dq^M$.
Then the following holds (for the statements ii. -- iv. we require that either $p+q\neq n-1,n$ or $p=0$):

\vspace{2mm}
i. $\Box^M_{p,q} = \dq_{p,q}^M (\dq_{p,q}^M)^* + (\dq_{p,q+1}^M)^* \dq_{p,q+1}^M$ has closed range
and the corresponding $\dq$-Neumann operator 
$$N^M_{p,q}=(\Box_{p,q}^M)^{-1}: L^{p,q}_\sigma(\Omega') \rightarrow L^{p,q}_\sigma(\Omega')$$
is bounded.\footnote{By a little abuse of notation,
we write $N=\Box^{-1}$ for $\dq$-Neumann operators though $N$ is an inverse to the $\dq$-Laplacian
$\Box=\dq\dq^*+\dq^*\dq$ only on the range $\Im \Box$.}

\vspace{2mm}
ii. On $\Omega^*=\Omega\setminus\Sing X$, the operators 
\begin{eqnarray*}
\dq_{p,q}: && L^{p,q-1}(\Omega^*) \rightarrow L^{p,q}(\Omega^*),\\
\dq_{p,q+1}: && L^{p,q}(\Omega^*) \rightarrow L^{p,q+1}(\Omega^*)
\end{eqnarray*}
both have closed range in the corresponding kernels of $\dq$.
$\Box_{p,q}=\dq_{p,q}\dq_{p,q}^* + \dq_{p,q+1}^*\dq_{p,q+1}$
has closed range and the corresponding $\dq$-Neumann operator
$$N_{p,q}=\Box_{p,q}^{-1}: L^{p,q}(\Omega^*) \rightarrow L^{p,q}(\Omega^*)$$ is bounded.

\vspace{2mm}
iii. If $\Im \dq_{p,q}^M$ is of finite codimension in $\ker\dq_{p,q+1}^M$,
then $\Im \dq_{p,q}$ is of finite codimension in $\ker \dq_{p,q+1}$ and both, $\Im \Box_{p,q}^M$ and $\Im \Box_{p,q}$,
are of finite codimension in $L^{p,q}_{\sigma}(\Omega')$ and $L^{p,q}(\Omega^*)$, respectively.
If $\Im \dq_{p,q+1}^M$ is of finite codimension in $\ker\dq_{p,q+2}^M$,
then $\Im \dq_{p,q+1}$ is of finite codimension in $\ker \dq_{p,q+2}$

\vspace{2mm}
iv. $N_{p,q}$ is compact on $\Omega^*$ exactly if $N^M_{p,q}$ is compact on $\Omega'$.
\end{thm}

We remark that closed range of $\dq^M_{p,q}$, $\dq^M_{p,q+1}$ in the assumptions
is in fact a necessary condition for (i) and (ii), respectively (see Theorem \ref{thm:fa6} and Theorem \ref{thm:compact}).

\smallskip
The present paper is organized as follows.
In Section \ref{sc:fa}, we generalize a characterization of precompact sets in the space
of $L^2$-forms on arbitrary Hermitian manifolds from \cite{Rp9} to the case of forms with values in
Hermitian line bundles (Theorem \ref{thm:precompact2})
which also induces a nice characterization of compactness of the $\dq$-Neumann operator
(Theorem \ref{thm:compact3}).
Besides, we collect some other facts from functional analysis that are needed throughout the paper.
In Section \ref{sec:resolution}, we recall the necessary facts about 
the resolution of singularities $\pi: M\rightarrow X$.

In Section \ref{sec:bundles}, we study the closed range property of the $\dq$-operator
and the existence of compact $\dq$-solution operators for forms with values in holomorphic line bundles
twisted along the exceptional set of the resolution $\pi: M\rightarrow X$ (Theorem \ref{thm:range1}).
We need this because $L^2$-forms on $X$ do not correspond to $L^2$-forms on the resolution $M$ in general,
but to $L^2$-forms with values in some line bundles (depending on the degree of the forms).

In Section \ref{sec:sop}, we recall the $\dq$-solution operators at isolated singularities
of Forn{\ae}ss, {\O}vrelid and Vassiliadou from \cite{FOV2},
and deduce in Theorem \ref{thm:sop} the existence of compact solution operators by use of the criterion 
given in Theorem \ref{thm:precompact2}.
These operators are used to prove our first main statement, Theorem \ref{thm:1.1}, in Section \ref{sec:global}.

Finally, in the last section, we use Theorem \ref{thm:1.1} and the representation
of the $\dq$-Neumann operator in terms of canonical $\dq$-solution operators (Theorem \ref{thm:N100})
to prove our second main statement, Theorem \ref{thm:1.2}.

\smallskip

%\newpage
\section{Preliminaries}\label{sc:fa}

\subsection{Compactness in the space of $L^2$-forms}\label{sec:compactness}

Let $M$ be a Hermitian manifold with volume form $dV_M$ and $(E,h)$ a Hermitian vector bundle over $M$.
If $f$ is a differential form on $M$ with values in $E$, we denote by $|f|_h$ its pointwise norm.
We denote by $L^{p,q}(M,E)$ the Hilbert space of locally measurable $(p,q)$-forms such that
$$\|f\|^2_{L^{p,q}(M,E)} := \int_M |f|^2_h dV_M < \infty.$$
For functions, we also write $L^2(M,E)$ instead of $L^{0,0}(M,E)$.

Consider a countable open covering $\{U_j\}_{j\in\N}$ of $M$. Then:

\begin{thm}\label{thm:precompact1}
A subset $\mathcal{K}$ of $L^2(M,E)$ is relatively compact exactly if
the following two conditions are fulfilled:

(i) $\mathcal{K}|_{U_j} =\{f|_{U_j}: f\in \mathcal{K}\}$ is relatively compact in $L^2(U_j,E|_{U_j})$ for every j.

(ii) For every $\epsilon>0$, there is $N\in\N$ such that
$$\int_{M\setminus\bigcup_{j=1}^N U_j} |f|^2_h dV_M < \epsilon\ \ \forall f\in \mathcal{K}.$$
\end{thm}

\begin{proof}
We use the elementary fact that a subset $\mathcal{K}$ of a metric space $\mathcal{M}$ is relatively compact
iff every sequence in $\mathcal{K}$ has a convergent subsequence (in $\mathcal{M}$).

\smallskip
Assume first that (i) and (ii) hold. Thus, let $\mathfrak{s}$ be a sequence in $\mathcal{K}$.
From (i), it follows by a Cantor diagonal construction that $\mathfrak{s}$ has a subsequence $\{f_\mu\}_\mu$
such that $\{f_\mu|_{U_j}\}_\mu$ converges in $L^2(U_j,E)$ for all $j$.

Let $\epsilon>0$. Choose $N$ so big that
$$\int_{M\setminus \bigcup_{j=1}^N U_j} |f|^2_h dV_M < \epsilon/8$$
for every $f\in\mathcal{K}$. Then
\begin{eqnarray*}
\int_M |f_{\mu} - f_{\nu}|^2_h dV_h 
&\leq& \sum_{j=1}^N \int_{U_j} |f_{\mu} - f_{\nu}|^2_h dV_M + 2 \int_{M\setminus \bigcup_{j=1}^N} \big(|f_\mu|^2 + |f_\nu|^2\big) dV_M\\
&\leq& \sum_{j=1}^N \int_{U_j} |f_{\mu} - f_{\nu}|^2_h dV_M + \epsilon/2,
\end{eqnarray*}
while the first term on the right hand side is less than $\epsilon/2$ when $\mu,\nu\geq \mu_0$ big enough.
So, $\{f_\mu\}_\mu$ is a Cauchy sequence in $L^2(M,E)$.

\smallskip

For the converse direction of the statement, assume that $\mathcal{K}$ is relatively compact.
Then property (i) is trivial. Let $\epsilon>0$. As $\overline{\mathcal{K}}$ is compact, there
exist finitely many functions $f_1, ..., f_K\in L^2(M,E)$ such that
\begin{eqnarray}\label{eq:coverballs}
\mathcal{K} \subset \bigcup_{\mu=1}^K B_{\sqrt{\epsilon}/2}(f_\mu)\ \ \mbox{ in } L^2(M,E).
\end{eqnarray}
We can now choose $N$ so big that
\begin{eqnarray*}
\int_{M\setminus \bigcup_{j=1}^N U_j} |f_\mu|^2_h dV_M < \epsilon/4 \ \ \mbox{ for }\mu=1, ..., K
\end{eqnarray*}
because we have to consider only finitely many functions simultaneously and can exhaust $M$ by increasing $N$ (see e.g. \cite{Alt}, A.1.16.2).
Let $f\in\mathcal{K}$. Then there exists by \eqref{eq:coverballs} an index $\mu_0$ such that
\begin{eqnarray*}
\int_{M\setminus \bigcup_{j=1}^N U_j} |f|^2_h dV_M
&\leq&
2 \int_{M\setminus \bigcup_{j=1}^N U_j} |f_{\mu_0}|^2_h dV_M + 
2 \int_{M\setminus \bigcup_{j=1}^N U_j} |f-f_{\mu_0}|^2_h dV_M\\
&<& \epsilon/2 + 2 \|f-f_{\mu_0}\|^2_{L^2(M,L)} < \epsilon,
\end{eqnarray*}
and that proves (ii).
\end{proof}

We can now use the G\r{a}rding inequality to characterize (relative) compactness of subsets of $L^2$-forms
which are bounded in the graph norm for $\dq\oplus \dq^*$.
To make that precise, we define
\begin{eqnarray*}
\|f\|^2_{\Gamma_{p,q}(M,E)} := \|f\|_{L^{p,q}(M,E)}^2 + \|\dq f\|^2_{L^{p,q+1}(M,E)} +\|\dq^* f\|^2_{L^{p,q-1}(M,E)}
\end{eqnarray*}
for $(p,q)$-forms $f \in L^{p,q}(M,L) \cap \Dom\dq\cap\Dom\dq^*$.
Here, $\dq$ is the $\dq$-operator in the sense of distributions 
$$\dq: L^{p,q}(M,E) \rightarrow L^{p,q+1}(M,E),$$
and 
$$\dq^*: L^{p,q}(M,E) \rightarrow L^{p,q-1}(M,E)$$
is the Hilbert space adjoint of $\dq: L^{p,q-1}(M,E)\rightarrow L^{p,q}(M,E)$.

Note that for forms of different degree with values in $E$,
we can use different Hermitian metrics for the line bundle $E$,
so that the spaces of $(p,q-1)$, $(p,q)$, $(p,q+1)$-forms carry different weights.

%The crucial characterization is:

\begin{thm}\label{thm:precompact2}
Let $M$ be a Hermitian manifold, and $E$ a Hermitian line bundle over $M$.
Let
$$\mathcal{A} \subset L^{p,q}(M,E)\cap \Dom\dq\cap\Dom\dq^*$$
be a set of $(p,q)$-forms which is bounded in the graph norm $\|\cdot\|_{\Gamma_{p,q}(M,E)}$.
Then $\mathcal{A}$ is relatively compact in $L^{p,q}(M,E)$ iff the following condition is fulfilled:

(C) For every $\epsilon>0$, there exists $\Omega_\epsilon\subset\subset M$ such that
$$\|f\|_{L^{p,q}(M-\Omega_\epsilon,E)}<\epsilon\ \ \mbox{ for all } f\in \mathcal{A}.$$
%for all $f\in\mathcal{A}$.
\end{thm}

\begin{proof}
The proof follows from Theorem \ref{thm:precompact1} by use of the G\r{a}rding inequality.
We can use Theorem \ref{thm:precompact1} for the bundle $\Lambda^{p,q}\otimes E$ so that $L^{p,q}(M,E) \cong L^2(M,\Lambda^{p,q}\otimes E)$.

We can choose a covering $\{U_j\}_j$ for $M$ such that the condition (i) in Theorem \ref{thm:precompact1}
is fulfilled for the set $\mathcal{A}$ because it is bounded in the graph norm.
Just make sure that $\{U_j\}_j$ is a countable covering such that each open set $U_j$ is relatively compact in $M$
and that $U_j$ is biholomorphic to a ball in $\C^n$. We can assume that $E$ is trivial over $U_j$.

It follows from the G\r{a}rding inequality (see e.g. \cite{FK}, Theorem 2.2.1) that there exists a constant $C_j>0$ such that
$$\|u\|_{W^{1,2}_{p,q}(U_j)} \leq C_j \|u\|_{\Gamma_{p,q}(M,E)}$$
for all $u\in C^\infty(M,E)$, where we denote by $W^{1,2}$ the Sobolev $W^{1,2}$-norm.
It follows by a standard density argument that $\mathcal{A}|_{U_j}$ is a bounded subset of $W^{1,2}_{p,q}(U_j)$.
But the embedding $W^{1,2}_{p,q}(U_j) \hookrightarrow L^{p,q}(U_j)$ is compact by the Rellich embedding theorem
(see \cite{Alt}, Theorem A.6.4).
Thus $\mathcal{A}|_{U_j}$ is relatively compact in $L^{p,q}(U_j,E|_{U_j})$ for each $j$.

The proof is completed by the easy observation that condition (C) is equivalent to condition (ii) from Theorem \ref{thm:precompact1}.
\end{proof}

Note that we intend to use Theorem \ref{thm:precompact2} with $M=\Reg X$, the regular set of a singular Hermitian space,
or an open subset of $\Reg X$.

From Theorem \ref{thm:precompact2}, one can deduce the following criterion for compactness of the $\dq$-Neumann operator
(see \cite{Rp9}, Theorem 1.3). Recall that the $\dq$-Neumann operator $N$ is defined as follows:
for $u\in\Im\Box$, let $Nu$ be the unique form in $\Box^{-1}(\{u\})$ which is orthogonal to $\ker\Box$.

\begin{thm}\label{thm:compact3}
Let $Z$ be a Hermitian complex space of pure dimension $n$, $X\subset Z$ an open Hermitian submanifold
and $\dq$ a closed $L^2$-extension of the $\dq_{cpt}$-operator on smooth forms with compact support in $X$,
for example the $\dq$-operator in the sense of distributions.
Let $0\leq p,q \leq n$.

Assume that $\dq$ has closed range in $L^{p,q}(X)$ and in $L^{p,q+1}(X)$.
Then $\Box = \dq \dq^* + \dq^*\dq$
has closed range in $L^{p,q}(X)$ and the following conditions are equivalent:

(i) The $\dq$-Neumann operator $N=\Box^{-1}: \Im\Box\rightarrow L^{p,q}(X)$ is compact.

(ii) For all $\epsilon>0$, there exists $\Omega\subset\subset X$ such that $\|u\|_{L^{p,q}(X-\Omega)}<\epsilon$
for all
$$u\in \{u\in \Dom(\dq)\cap\Dom(\dq^*)\cap \Im \Box: \|\dq u\|_{L^{p,q+1}} + \|\dq^*u\|_{L^{p,q-1}}<1\}.$$

(iii) There exists  a smooth function $\psi\in C^\infty(X,\R)$, $\psi>0$, such that $\psi(z)\rightarrow \infty$ as $z\rightarrow bX$,
and
\begin{eqnarray*}
( \Box u,u)_{L^2} \geq \int_X \psi |u|^2 dV_X\ \ \mbox{ for all } u\in\Dom\Box\cap\Im\Box \subset L^{p,q}(X).
\end{eqnarray*}
\end{thm}

\medskip

\subsection{Closed, densely defined linear operators}
\label{ssec:linops}

Let $T: H_1 \rightarrow H_2$
be a closed, densely defined linear operator between Hilbert spaces $H_1$, $H_2$.
Then its adjoint operator $T^*: H_2 \rightarrow H_1$ has the same properties.

We define the {\bf minimal solution operator} of $T$,
${\bf S}: \Im T \rightarrow H_1$,
by $T\circ {\bf S} (u)=u$ and ${\bf S}u\perp \ker T$, i.e. ${\bf S}u$ is the unique
element of $T^{-1}(\{u\})$ that is perpendicular to $\ker T$.
It follows directly by \cite{Hoe1}, Theorem 1.1.1, that ${\bf S}$ is bounded exactly if
$\Im T$ is closed. Moreover, $\Im T$ is closed exactly if $\Im T^*$ is closed.

When $\Im T$ is closed, we extend the minimal operator ${\bf S}$ to a bounded operator $H_2 \rightarrow H_1$
by setting ${\bf S}u=0$ for $u\in (\Im T)^\perp=\ker T^*$.

\begin{lem}\label{lem:fa4}
Assume that $\Im T$ is closed and let ${\bf S}: H_2 \rightarrow H_1$ be the minimal solution operator of $T$.
Then the minimal solution operator of $T^*$ equals the adjoint ${S}^*: H_1 \rightarrow H_2$ of ${\bf S}$.
\end{lem}

\begin{proof}
When ${\bf S}'$ is the minimal solution operator of $T^*$, we observe
for $u\in H_1$ and $v\in H_2$ that
$${\bf S}v \in (\ker T)^\perp\ \ \mbox{ and } \ \ {\bf S}' u\in (\ker T^*)^\perp.$$
Moreover,
$$v-T{\bf S} v \in (\Im T)^\perp=\ker T^* \ \ \mbox{ and } \ \ u-T^*{\bf S}'u \in (\Im T^*)^\perp=\ker T.$$
Hence
$$({\bf S}'u,v)_{H_2} = ({\bf S}'u,T{\bf S}v)_{H_2} = (T^* {\bf S}'u,{\bf S}v)_{H_1} = (u,{\bf S}v)_{H_1},$$
since $v-T{\bf S}v \perp {\bf S}'u$ and $u-T^*{\bf S}' \perp {\bf S}v$. Thus ${\bf S}'={\bf S}^*$.
\end{proof}

Consider next two closed, densely defined linear Hilbert space operators
$T_1: H_1 \rightarrow H_2$, $T_2: H_2 \rightarrow H_3$
such that $\Im T_1 \subset \ker T_2$, i.e. $T_2\circ T_1=0$.
Then we have also that $\Im T_2^*\subset \ker T_1^*$, i.e. $T_1^*\circ T_2^*=0$.
We define
$$\Box: H_2\rightarrow H_2 \ \ \mbox{ by } \ \ \Box u=T_1T_1^* u + T_2^* T_2 u$$
on $\Dom \Box =\{ u\in \Dom T_1^*\cap \Dom T_2: T_1^* u \in\Dom T_1\ \mbox{ and }\ T_2u\in \Dom T_2^*\}$.

\begin{thm}\label{thm:fa5}
$\Box$ is closed, densely defined and self-adjoint.
\end{thm}

The statement is well-known fact from operator theory, for a proof see e.g. \cite{Rp9}, Theorem 3.1.
Clearly, the kernel of $\Box$ is 
$$\mathcal{H}=\ker T_1^*\cap \ker T_2$$
since $(\Box u,u)=\|T_1^* u\|^2 + \|T_2 u\|^2$ when $u\in \Dom\Box$, and we have
the orthogonal decomposition
$$H_2 = \mathcal{H} \oplus^\perp \o{\Im \Box}$$
by the self-adjointness of $\Box$. Observe that $\mathcal{H}$, $\Im T_1$ and $\Im T_2^*$ are mutually orthogonal.
We have the obvious inclusions
\begin{eqnarray}\label{eq:oi}
\o{\Im\Box} \supset \Im T_1 + \Im T_2^* \supset \Im T_1T_1^* + \Im T_2^*T_2 \supset \Im \Box.
\end{eqnarray}
Taking closures gives the orthogonal decompositions
\begin{eqnarray*}
\o{\Im \Box} = \o{\Im T_1} \oplus^\perp \o{\Im T_2^*} = \o{\Im T_1T_1^*}\oplus^\perp \o{\Im T_2^*T_2},
\end{eqnarray*}
and
\begin{eqnarray*}
H_2 = \mathcal{H}\oplus^\perp \o{\Im T_1T_1^*} \oplus^\perp \o{\Im T_2^* T_2}.
\end{eqnarray*}

\begin{thm}\label{thm:fa6}
$\Box$ has closed image exactly if $T_1$ and $T_2$ have closed images.
Let $N$ be the minimal solution operator of $\Box$.
When $\Im \Box$ is closed, then
$$N= {\bf S}^*_1 {\bf S}_1 + {\bf S}_2 {\bf S}^*_2$$
for the extension of $N$ to $H_2$.
\end{thm}

\begin{proof}
If $\Im \Box$ is closed, then \eqref{eq:oi} implies that
$\Im \Box=\Im T_1\oplus^\perp\Im T_2^*$, where $\Im T_1$ and $\Im T_2^*$ must be closed.

On the other hand, assume that $\Im T_1$ and $\Im T_2^*$ are closed.
Let $f=f_1 + f_2$ where $f_1\in \Im T_1$ and $f_2\in \Im T_2^*$.
Then
$${\bf S}_1 f = {\bf S}_1 f_1 \in (\ker T_1)^\perp=\Im T_1^*$$
since $\Im T_1^*$ is closed. Also
$$u_1= {\bf S}_1^* {\bf S}_1 f_1 = {\bf S}_1^* {\bf S}_1 f \in (\ker T_1^*)^\perp = \Im T_1$$
so that $T_2 u_1=0$ while $T_1T_1^* u_1=f_1$. In the same way, 
$$u_2={\bf S}_2{\bf S}_2^* f$$
satisfies $T_2^* T_2 u_2=f_2$ and $T_1^* u_2=0$.
Setting $u=u_1+u_2$, we see that 
$$u\in\Dom\Box\ \ \mbox{ and }\ \ \Box u=f$$
so that $f\in\Im\Box$. On the other hand, $\Im \Box\subset \Im T_1 + \Im T_2^*$,
so that the both spaces must be equal.

Moreover, $u_1\in\Im T_1$ and $u_2\in \Im T_2^*$, so $u\perp \mathcal{H}$
and thus 
$$u=({\bf S}_1^*{\bf S}_1 + {\bf S}_2 {\bf S}_2^*) f = N f.$$
\end{proof}

The first part of the proof of the Theorem implies the following:

\begin{prop}\label{prop:fa7}
If $T_1$ has closed image, then 
$$\Im\Box \supset \Im T_1\ \ \mbox{ and }\ \ N|_{\Im T_1}={\bf S}_1^* {\bf S}_1 |_{\Im T_1}.$$
\end{prop}

\medskip

\subsection{Some more functional analysis}\label{sec:fa}
In this section, we state some well-known results from functional analysis 
that are used in the paper. We include proofs for convenience of the reader when we do not know a suitable reference.

\begin{lem}\label{lem:fa1}
If $V$ is a subspace of a Banach space $B$ that has a closed subspace $V_0$ of finite codimension in $V$,
then $V$ is also closed.
\end{lem}

\begin{proof}
As $V_0$ is closed, the quotient space $B/V_0$ is also a Banach space and the quotient map $q:  B \rightarrow B/V_0$ is continuous
(see e.g. \cite{Rd}, Theorem 1.41).

$B/V$ is a finite dimensional subspace of $B/V_0$, hence closed.
So, $V=q^{-1}(B/V)$ is also closed.
\end{proof}

\begin{lem}\label{lem:fa2}
Let $A: E\rightarrow F$ be a closed operator between Banach spaces $E$ and $F$ with the domain of definition $\Dom A\subset E$.
If $A(\Dom A)$ is of finite codimension in $F$, then $A(\Dom(A))$ is a closed subspace of $F$.
\end{lem}

Lemma \ref{lem:fa2} follows from Banach's open mapping theorem. For a proof, see \cite{HeLe}, Appendix 2.4.

\begin{lem}\label{lem:fa3}
Let $T: H_1 \rightarrow H_2$ be a linear map between Hilbert spaces $H_1, H_2$.
If $T|_V: V \rightarrow H_2$ is compact for a closed subspace $V$ of finite codimension in $H_1$,
then $T$ is compact.
\end{lem}

\begin{proof}
Let $\{a_1, ..., a_k\}$ be a basis for the orthogonal complement of $V$ in $H_1$,
and let $b_\mu:= T a_\mu$. Let $\{f_j\}_j$ be a bounded sequence in $H_1$ and consider the unique
representation
$$f_j = f_j^0 + f_j^1 a_1 + \cdots +f_j^k a_k$$
where $f_j^0 \in V$ and $f_j^1, ..., f_j^k \in \C$. As $T|_V$ is compact, there exists a subsequence $\{f_{j_l}\}_l$ of $\{f_j\}_j$ such that $\{Tf_{j_l}^0\}_l$
is a Cauchy sequence in $H_2$. 

On the other hand $\{f_j\}_j$ is bounded so that $\{f_j^\mu\}_j$ is bounded for $\mu=1, ..., k$.
So, we can choose the subsequence $\{f_{j_l}\}_l$ such all that the sequences $\{f_{j_l}^\mu\}_l$
are convergent for $\mu=1, ..., k$. But then
$$Tf_{j_l} = Tf_{j_l}^0 + f_{j_l}^1 b_1 + \cdots + f_{j_l}^k b_k$$
is also a Cauchy sequence in $H_2$.
\end{proof}

\begin{lem}\label{lem:fa5}
Let $T: H_1 \rightarrow H_2$ be a closed and surjective linear map between Hilbert spaces $H_1, H_2$,
and ${\bf S}_0: V \rightarrow H_1$ a bounded right-inverse to $T$ on a subspace $V$ of finite codimension in $H_2$.
Then ${\bf S}_0$ can be extended to a right-inverse $H_2\rightarrow H_1$ to $T$, and any such extension is bounded.
If ${\bf S}_0$ is compact, then so is ${\bf S}$.
\end{lem}

\begin{proof}
Choose a basis $e_1, ..., e_l$ of the complement of $V$ in $H_2$.
Then there exists forms $h_1, ..., h_l\in H_1$
such that $T h_\nu =e_\nu$ for $\nu=1, ..., l$.
Then each $f\in H_2$ has a unique representation
\begin{eqnarray*}
f=f' + \sum_{\nu=1}^l a_\nu e_\nu,\ \ \ f'\in V,\  a_1, ..., a_l\in\C,
\end{eqnarray*}
and we define ${\bf S} (f) := {\bf S}_0(f') + \sum_{\nu=1}^l a_\nu h_\nu$.
It is clear that any extension of ${\bf S}_0$ to $H_2$ has such a representation and is bounded.
If ${\bf S}_0$ is compact, then ${\bf S}$ is compact by Lemma \ref{lem:fa3}.
\end{proof}

\smallskip

\subsection{Resolution of singularities and comparison of metrics}\label{sec:resolution}

Let $\pi: M \rightarrow X$
be a resolution of singularities (which exists due to Hironaka \cite{Hi}), i.e. a proper holomorphic surjection such that
$\pi|_{M-E}: M-E \rightarrow X-\Sing X$
is biholomorphic, where $E=\pi^{-1}(\Sing X)$ is the exceptional set.
We may assume that $E$ is a divisor with only normal crossings,
i.e. the irreducible components of $E$ are regular and meet complex transversely.
However, this assumption is not really necessary for the results of this paper,
it is enough to assume that $E$ is a divisor.
For the topic of desingularization, we refer to
\cite{AHL}, \cite{BiMi} and \cite{Ha}.
Let $\gamma:= \pi^* h$
be the pullback of the Hermitian metric $h$ of $X$ to $M$.
$\gamma$ is positive semidefinite (a pseudo-metric) with degeneracy locus $E$.

We give $M$ the structure of a Hermitian manifold with a freely chosen (positive definite)
metric $\sigma$. Then $\gamma \lesssim \sigma$
and $\gamma \sim \sigma$ on compact subsets of $M-E$.
For an open set $U\subset M$, we denote by $L^{p,q}_{\gamma}(U)$ and $L^{p,q}_{\sigma}(U)$
the spaces of square-integrable $(p,q)$-forms with respect to the (pseudo-)metrics $\gamma$ and $\sigma$,
respectively.

For an open set $\Omega \subset X$, $\Omega^*=\Omega - \Sing X$, $\Omega':=\pi^{-1}(\Omega)$,
pullback of forms under $\pi$ gives the isometry
\begin{eqnarray}\label{eq:l2est5}
\pi^*: L^{p,q}(\Omega^*) \rightarrow L^{p,q}_{\gamma}(\Omega'-E) \cong L^{p,q}_{\gamma}(\Omega'),
\end{eqnarray}
where the last identification is by trivial extension of forms over the thin exceptional set $E$.

\newpage
\section{Holomorphic line bundles twisted along the exceptional set}\label{sec:bundles}

Let $D$ be a divisor on $M$ with support on the exceptional set $E$ of the resolution $\pi: M\rightarrow X$,
and denote by $\OO(D)$ the sheaf of germs of meromorphic functions $f$ such that $\mbox{div}(f)+D\geq 0$.
We denote by $L_D$ the associated holomorphic line bundle such that sections in $\OO(D)$ correspond
to holomorphic sections in $L_D$.
The constant function $f\equiv 1$ induces a meromorphic section $s_D$ of $L_D$
such that $\mbox{div}(s_D)=D$. One can then identify sections in $\OO(D)$
with sections in $\OO(L_D)$ by $g\mapsto g\otimes s_D$,
and we denote the inverse mapping by $g\mapsto g\cdot s_D^{-1}$.
If $D$ is an effective divisor, then $s_D$ is a holomorphic section of $L_D$ and $\OO(-D) \subset \OO \subset \OO(D)$.
If $Z$ is an effective divisor, then there is the natural inclusion $\OO(D)\subset \OO(D+Z)$
which induces the injection $\OO(L_{D})\subset \OO(L_{D+Z})$ given by
$g\mapsto (g\cdot s_{D}^{-1})\otimes s_{D+Z}$.
This also induces the natural injection of smooth sections of vector bundles
\begin{eqnarray}\label{eq:inclusion1}
\Gamma(U,L_{D}) \subset \Gamma(U,L_{D+Z})
\end{eqnarray}
for open sets $U\subset M$.

We give each $L_D$ the structure of a Hermitian holomorphic line bundle by choosing an arbitrary
positive definite Hermitian metric on $L_D$. We denote by
$$L^{p,q}_{\sigma}(U,L_D)$$
the space of square-integrable $(p,q)$-forms with values in $L_D$ (with respect to the metric $\sigma$
on $M$ and the chosen metric on $L_D$).
If $U\subset\subset M$ is relatively compact, then \eqref{eq:inclusion1}
induces the natural injection
\begin{eqnarray}\label{eq:inclusion2}
L^{p,q}_{\sigma}(U,L_D) \subset L^{p,q}_{\sigma}(U,L_{D+Z})
\end{eqnarray}
for any effective divisor $Z$. This does not depend on the metrics chosen on the line bundles
$L_D$ and $L_{D+Z}$ because $U$ is relatively compact in $M$.
For more details on Hermitian holomorphic
line bundles twisted along the exceptional set of the desingularization, we refer to Section 2 in \cite{Rp8} and Section 2 in \cite{OV3}.

%Note that the inclusions \eqref{eq:l2est3} and \eqref{eq:l2est4} remain valid for forms with values
%in a line bundle $L_D$. Moreover,

By \cite{Rp1}, Lemma 2.1, or \cite{FOV1}, Lemma 3.1, respectively,
there exists an effective divisor $\wt{D}$ with support on the exceptional set $E$
such that
\begin{eqnarray}\label{eq:inclusion3}
L^{p,q}_{\sigma}(U,L_{-\wt{D}}) \subset L^{p,q}_{\gamma}(U) \subset L^{p,q}_\sigma(U,L_{\wt{D}})
\end{eqnarray}
for all $0\leq p,q\leq n$ and open sets $U\subset\subset M$. 
This follows from the fact that $dV_\gamma$ vanishes
of a certain order (exactly) on $E$. 
For simplicity, we assume that $X$ has only
finitely many isolated singularities so that we can choose a fixed positive integer
$m$ such that the effective divisor $mE$ satisfies \eqref{eq:inclusion3}:
\begin{eqnarray}\label{eq:inclusion4}
L^{p,q}_{\sigma}(U,L_{-mE}) \subset L^{p,q}_{\gamma}(U) \subset L^{p,q}_\sigma(U,L_{mE})
\end{eqnarray}
for all $0\leq p,q\leq n$ and open sets $U\subset\subset M$.

In the present paper, we are interested in properties of the $\dq$-operator
$$\dq: L^{p,q}_{\gamma}(U-E) \rightarrow L^{p,q+1}_{\gamma}(U-E),$$
which we would like to relate to properties of the $\dq$-operator
$$\dq: L^{p,q}_{\sigma}(U) \rightarrow L^{p,q+1}_{\sigma}(U),$$
where $U$ is a neighborhood of the exceptional set.
We denote by $\dq_{D}$ the $\dq$-operator acting (in the sense of distributions) on $L^2$-forms with values in $L_D$.

As a preparation, we need:

\begin{thm}\label{thm:range1}
Let $U\subset\subset M$ be a neighborhood of the exceptional set $E$ and $D_1$, $D_2$ two divisors
with support on $E$. Then
\begin{eqnarray}\label{eq:dq1}
\dq_{D_1}: L^{p,q}_{\sigma}(U,L_{D_1}) \rightarrow L^{p,q+1}_{\sigma}(U,L_{D_1})
\end{eqnarray}
has closed range (of finite codimension) in $\ker\dq_{D_1}$ exactly if
\begin{eqnarray}\label{eq:dq2}
\dq_{D_2}: L^{p,q}_{\sigma}(U,L_{D_2}) \rightarrow L^{p,q+1}_{\sigma}(U,L_{D_2})
\end{eqnarray}
has closed range (of finite codimension) in $\ker \dq_{D_2}$. If this is the case, then there exists a compact
$\dq$-solution operator 
\begin{eqnarray}\label{eq:dq11}
{\bf S}_1: \Im \dq_{D_1} \subset L^{p,q+1}_{\sigma}(U,L_{D_1}) \rightarrow L^{p,q}_{\sigma}(U,L_{D_1})
\end{eqnarray}
exactly if there exists a compact $\dq$-solution operator
\begin{eqnarray}\label{eq:dq12}
{\bf S}_2: \Im \dq_{D_2} \subset L^{p,q+1}_{\sigma}(U,L_{D_2}) \rightarrow L^{p,q}_{\sigma}(U,L_{D_2}).
\end{eqnarray}
\end{thm}

\begin{proof}
It is enough to prove the statement for two divisors $D_1$, $D_2$ such that $D_2 - D_1\geq 0$.
The general statement follows then by comparing both bundles, $L_{D_1}$ and $L_{D_2}$, with $L_{|D_1|+|D_2|}$,
because $(|D_1|+|D_2|) -D_j\geq 0$ for $j=1,2$.

Since $E$ is compact by assumption (because $E\subset U\subset\subset M$),
it consists of $k$ pairwise disjoint components $E_\mu$, $\mu=1, ..., k$,
such that each $a_\mu=\pi(E_\mu)$ is an isolated singularity. 
For $\mu=1, ..., k$ choose pairwise disjoint strongly pseudoconvex neighborhoods
$V_\mu\subset\subset U$ of the components $E_\mu$. Let $V:=\bigcup_\mu V_\mu$.

Then it is well known that the operators
\begin{eqnarray*}
\dq_{D_1}|_{V}: L^{p,q}_{\sigma}(V,L_{D_1}) &\rightarrow& L^{p,q+1}_{\sigma}(V,L_{D_1}),\\
\dq_{D_2}|_{V}: L^{p,q}_{\sigma}(V,L_{D_2}) &\rightarrow& L^{p,q+1}_{\sigma}(V,L_{D_2})
\end{eqnarray*}
have closed range of finite codimension in the corresponding kernels of $\dq_{D_1}|_{V}$ and $\dq_{D_2}|_{V}$,
respectively, and that there are corresponding compact $\dq$-solution operators
(see e.g. the proof of Lemma 2.2 in \cite{OV3} which implies that the corresponding $\dq$-Neumann
operators are compact). 
We denote by $H_1^V$ and $H_2^V$ the range of $\dq_{D_1}|_{V}$ and $\dq_{D_2}|_{V}$
in $L^{p,q+1}_{\sigma}(V,L_{D_1})$ and $L^{p,q+1}_{\sigma}(V,L_{D_2})$, respectively,
and by
\begin{eqnarray*}
{\bf T}_1: H_1^V \rightarrow L^{p,q}_{\sigma}(V,L_{D_1})\ ,\ \
{\bf T}_2: H_2^V \rightarrow L^{p,q}_{\sigma}(V,L_{D_2})
\end{eqnarray*}
corresponding compact $\dq$-solution operators.

\medskip
Assume first that $\dq_{D_1}$ has closed range of finite codimension
and that
\begin{eqnarray}\label{eq:solution1}
{\bf S}_1: \Im \dq_{D_1} \subset L^{p,q+1}_{\sigma}(U,L_{D_1}) \rightarrow L^{p,q}_{\sigma}(U,L_{D_1})
\end{eqnarray}
is a corresponding bounded $\dq$-solution operator.
Consider the bounded linear map
$$\Phi_2: \ker \dq_{D_2} \subset L^{p,q+1}_{\sigma}(U,L_{D_2}) \rightarrow \ker \dq_{D_2}|_{V}$$
given by $\Phi_2(f) := f|_{V}$.
Since $H_2^V=\Im \dq_{D_2}|_V$ is a closed subspace of finite codimension in $ \ker \dq_{D_2}|_{V}$,
the same holds for
$$H_2:=\Phi_2^{-1}\big(H^V_2\big) = \big\{ f\in \ker \dq_{D_2} : f|_{V} \in H^V_2 \big\}$$
in $\ker \dq_{D_2}$. Now then, choose a smooth cut-off function $\chi$ which has compact support in $V=\bigcup_\mu V_\mu$
and is identically $1$ in a smaller neighborhood of the exceptional set $E$.
Then we can define a bounded linear map
$$\Psi_2: H_2 \rightarrow \ker \dq_{D_1} \subset L^{p,q+1}_{\sigma}(U,L_{D_1})$$
by the assignment
$$f\in H_2 \mapsto \big(\big(f - \dq_{D_2} ( \chi {\bf T}_2 ( f|_{V}))\big)\cdot s_{D_2}^{-1} \big)\otimes s_{D_1}$$
since $f - \dq_{D_2} ( \chi {\bf T}_2 ( f|_{V}))$
is identically zero in a neighborhood of $E$ and $\supp D_1\subset E$.
By assumption, $\dq_{D_1}$ has closed range $\Im \dq_{D_1}$ of finite codimension in $\ker\dq_{D_1}$
so that
$$H_2':=\Psi_2^{-1} ( \Im \dq_{D_1})$$
has closed range of finite codimension in $H_2$. 
As we have already seen that $H_2$ in turn is closed of finite codimension in $\ker\dq_{D_2}$,
it follows that $H_2'$ is a closed subspace of finite codimension in $\ker\dq_{D_2}$.

On the other hand, since $\Psi_2(H_2') \subset \Im\dq_{D_1}$,
we can define by use of \eqref{eq:solution1} a $\dq$-solution operator
\begin{eqnarray*}
{\bf S}_2': H_2' \subset L^{p,q+1}_{\sigma}(U,L_{D_2}) \rightarrow L^{p,q}_{\sigma}(U,L_{D_2})
\end{eqnarray*}
by setting
\begin{eqnarray*}
{\bf S}_2' (f):= \big(({\bf S_1}\circ\Psi_2(f))\cdot s_{D_1}^{-1}\big)\otimes s_{D_2} + \chi {\bf T}_2(f|_{V}).
\end{eqnarray*}
Here, we use the natural injection $\Gamma(U,L_{D_1}) \subset \Gamma(U,L_{D_2})$, $g\mapsto (g\cdot s_{D_1}^{-1})\otimes s_{D_2}$,
which exists since $D_2\geq D_1$ by assumption.
Since this injection commutes with the $\dq$-operator,
\begin{eqnarray*}
\dq_{D_2} \big(  \big(({\bf S}_1\circ\Psi_2(f))\cdot s_{D_1}^{-1}\big)\otimes s_{D_2}\big) 
&=& \big(\big(\dq_{D_1} ({\bf S}_1\circ\Psi_2(f))\big)\cdot s_{D_1}^{-1}\big)\otimes s_{D_2}\\
&=& \big( \Psi_2(f) \cdot s_{D_1}^{-1}\big)\otimes s_{D_2}\\
&=& f - \dq_{D_2} ( \chi {\bf T}_2 ( f|_{V})).
\end{eqnarray*}
Hence, $\dq_{D_2} {\bf S}_2'(f)=f$, so that $H_2'\subset \Im\dq_{D_2}$. Summing up, $H_2' \subset \Im\dq_{D_2} \subset \ker \dq_{D_2}$,
where $H_2'$ is closed and of finite codimension in $\ker \dq_{D_2}$. 
But then $\Im\dq_{D_2}$ is also closed and of finite codimension in $\ker\dq_{D_2}$ (see Lemma \ref{lem:fa1}).
By Lemma \ref{lem:fa5}, we can now extend ${\bf S}_2'$ to a $\dq$-solution operator ${\bf S}_2$ on $\Im\dq_{D_2}$.

If, in addition, ${\bf S}_1$ is compact, then ${\bf S}_2'$ is compact
because ${\bf T}_2$ is also compact.
But then ${\bf S}_2$ is a compact $\dq$-solution operator on $\Im \dq_{D_2}$ (see Lemma \ref{lem:fa3}).

\medskip

For the converse direction of the statement, assume that 
that $\dq_{D_2}$ has closed range of finite codimension
and that
\begin{eqnarray}\label{eq:solution2}
{\bf S}_2: \Im \dq_{D_2} \subset L^{p,q+1}_{\sigma}(U,L_{D_2}) \rightarrow L^{p,q}_{\sigma}(U,L_{D_2})
\end{eqnarray}
is a corresponding bounded $\dq$-solution operator.
We have to distinguish the two cases $q\geq 1$ and $q=0$.

Assume first that $q\geq 1$.
Analogously to the converse direction above,
consider now the bounded linear map
$$\Phi_1: \ker \dq_{D_1} \subset L^{p,q+1}_{\sigma}(U,L_{D_1}) \rightarrow \ker \dq_{D_1}|_{V}$$
given by $\Phi_1(f) := f|_{V}$.
Since $H^V_1=\Im\dq_{D_1}|_V$ is a closed subspace of finite codimension in $ \ker \dq_{D_1}|_{V}$,
the same holds for
$$H_1:=\Phi_1^{-1}\big(H^V_1\big) = \big\{ f\in \ker \dq_{D_1} : f|_{V} \in H^V_1\big\}$$
in $\ker \dq_{D_1}$.
Again, we use the cut-off function $\chi$ to define a bounded linear map
$$\Psi_1: H_1 \rightarrow \ker \dq_{D_2} \subset L^{p,q+1}_{\sigma}(U,L_{D_2})$$
by the assignment
$$f\in H_1 \mapsto \big(\big(f - \dq_{D_1} ( \chi {\bf T}_1 ( f|_{V}))\big)\cdot s_{D_1}^{-1} \big)\otimes s_{D_2}.$$
Here now, it seems at first hand superfluous that
$f - \dq_{D_1} ( \chi {\bf T}_1 ( f|_{V}))$
is identically zero in a neighborhood of $E$ because we can use the natural injection
$\Gamma(U,L_{D_1}) \subset \Gamma(U,L_{D_2})$, $g\mapsto (g\cdot s_{D_1}^{-1})\otimes s_{D_2}$ (recall that $D_2\geq D_1$).
But, we have now to face the problem that forms in the image of ${\bf S}_2$ can not be considered
as forms with values in $L_{D_1}$. To treat this difficulty, let $W_1, ..., W_k$
be strongly pseudoconvex neighborhoods of the components $E_1, ..., E_k$ of the exceptional set
which are contained in the set where $\chi\equiv 1$, and let $W:=\bigcup_\mu W_\mu$.
Denote by $B_2^W$ the range of $\dq_{D_2}|_{W}$ in $L^{p,q}_{\sigma}(W,L_{D_2})$,
and let
$${\bf K}_2: B_2^W \rightarrow L^{p,q-1}_{\sigma}(W,L_{D_2})$$
be a corresponding compact $\dq$-solution operator. 

%Let
%$$\Phi_2': \ker \dq_{D_2}|_{W} \subset L^{p,q}_{\sigma}(W,L_{D_2}) \rightarrow \ker\dq_{D_2}|_{W}$$
%be the bounded linear map given by
%$$\Phi_2'(f) := f|_{W}.$$
%Since $B^W_2$ is a closed subspace of finite codimension in $G_2'$,
%$$B_2:=(\Phi_2')^{-1}\big(B^1_2\oplus \cdots \oplus B^k_2\big) = \big\{ f\in \ker \dq_{D_2}|_W : f|_{W_\mu} \in B^\mu_2, \mu=1, ..., k\big\}$$
%is a closed subspace of finite codimension in $\ker \dq_{D_2}|_{W} \subset L^{p,q}_{\sigma}(W,L_{D_2})$.

By assumption, $\dq_{D_2}$ has closed range $\Im \dq_{D_2}\subset L^{p,q+1}_{\sigma}(U,L_{D_2})$ of finite codimension in 
$\ker\dq_{D_2}\subset L^{p,q+1}_{\sigma}(U,L_{D_2})$
so that
$$H_1':=\Psi_1^{-1} ( \Im \dq_{D_2})$$
has closed range of finite codimension in $H_1$. 
As we have already seen that $H_1$ in turn is closed of finite codimension in $\ker\dq_{D_1}$,
it follows that $H_1'$ is a closed subspace of finite codimension in $\ker\dq_{D_1}$.

Consider
$${\bf S}_2 \circ \Psi_1|_{H_1'}: H_1' \rightarrow L^{p,q}_{\sigma}(U,L_{D_2}).$$
Since
$$\Psi_1(f) = \big(\big(f - \dq_{D_1} ( \chi {\bf T}_1 ( f|_{V}))\big)\cdot s_{D_1}^{-1} \big)\otimes s_{D_2}$$
vanishes on $W$, it follows that ${\bf S}_2\circ \Psi_1|_{H_1'}(f)$
is $\dq$-closed on $W$, and so ${\bf S}_2 \circ \Psi_1|_{H_1'}: H_1' \rightarrow \ker\dq_{D_2}|_W$
is a bounded linear map. Let
$$H_1'':= \big({\bf S}_2 \circ \Psi_1|_{H_1'}\big)^{-1}(B_2^W).$$
It follows that $H_1''$ is closed of finite codimension in $H_1'$ as $B_2^W$ is closed 
and of finite codimension in $\ker\dq_{D_2}|_W$. So, $H_1''$ is closed of finite codimension
in $\ker\dq_{D_1}$. 

On the other hand, $H_1''\subset \Im \dq_{D_1}$, as can be seen now as follows.
Let $\chi'$ be a smooth cut-off function with support in $W$ that is identically $1$
in a smaller neighborhood of the exceptional set $E$, and let
$${\bf S}_1'': H_1'' \rightarrow L^{p,q}_{\sigma}(U,L_{D_2})$$
be given as
\begin{eqnarray}\label{eq:S1''}
{\bf S}_1''(f):= {\bf S}_2 \circ \Psi_1|_{H_1'} (f) - \dq_{D_2} \chi' {\bf K}_2\big( \left.{\bf S}_2 \circ \Psi_1|_{H_1'} (f)\right|_{W}\big).
\end{eqnarray}
Since the forms ${\bf S}_1''(f)$ vanish identically in a neighborhood of the exceptional set $E$,
we can define
\begin{eqnarray}\label{eq:complement3}
{\bf S}_1': H_1'' \rightarrow L^{p,q}_{\sigma}(U,L_{D_1}),\ f\mapsto \big({\bf S}_1''(f)\cdot s_{D_2}^{-1}\big)\otimes s_{D_1},
\end{eqnarray}
and it is easy to check that $\dq_{D_1} {\bf S}_1'(f)=f$ for all $f\in H_1''$.
So, also $\Im \dq_{D_1}$ is closed and of finite codimension in $\ker\dq_{D_1}$
and one can complement ${\bf S}_1'$ by Lemma \ref{lem:fa5} to a bounded $\dq$-solution operator
\begin{eqnarray}\label{eq:complement4}
{\bf S}_1: \Im \dq_{D_1} \subset L^{p,q+1}_{\sigma}(U,L_{D_1}) \rightarrow L^{p,q}_{\sigma}(U,L_{D_1}).
\end{eqnarray}
If in addition ${\bf S}_2$ is compact, then ${\bf S}_1''$ and ${\bf S}_1'$ are compact, as well.
The additional fact that
${\bf T}_1$ and ${\bf K}_2$ are compact is not needed here since
the composition of a compact and a bounded operator is compact.
Note that here the right-hand side of \eqref{eq:S1''} does not cause any difficulties
because
\begin{eqnarray*}
 \dq_{D_2} \chi' {\bf K}_2\big( \left.{\bf S}_2 \circ \Psi_1|_{H_1'} (f)\right|_{W}\big) =
 \dq \chi' \wedge {\bf K}_2\big( \left.{\bf S}_2 \circ \Psi_1|_{H_1'} (f)\right|_{W}\big)  +\  \chi' {\bf S}_2 \circ \Psi_1|_{H_1'} (f).
\end{eqnarray*}
Consequently, ${\bf S}_1$ is a compact $\dq$-solution operator on $\Im \dq_{D_1}$.

\medskip

It remains to treat the case $q=0$. Let $\Phi_1, H_1, \Psi_1$ and $H_1'$ be as above (in the case $q\geq 1$).
We shall simply change the definition of $H_1''$:

Observe that when $f\in H_1'$, $\left.{\bf S}_2\circ\Psi_1(f)\right|_W$ is a holomorphic section of $\Omega^p\otimes L_{D_2}$ over $W$.
Let $\omega^p_{D_j}$ be the sheaf of holomorphic sections of $\Omega^p\otimes L_{D_j}$, $j=1,2$.
Then there is a short exact sequence of coherent analytic sheaves
$$0\rightarrow \omega^p_{D_1} \overset{i}{\longrightarrow} \omega^p_{D_2} \longrightarrow Q \rightarrow 0,$$
where $i(h)=(h\cdot s^{-1}_{D_1})\otimes s_{D_2}$ and $Q$ is supported on the exceptional set $E$.
This gives the exact sequence
$$0\rightarrow \Gamma(W,\omega^p_{D_1}) \overset{i_*}{\longrightarrow} \Gamma(W,\omega^p_{D_2}) \longrightarrow \Gamma(W,Q),$$
where $\Gamma(W,Q)$ is of finite dimension as $Q$ is a coherent analytic sheaf with compact support in $W$.
So, $\Im i_*$ has finite codimension in $\Gamma(W,\omega^p_{D_2})$.

We have
$$\ker \dq^{p,1}_{D_j}(W)= \Gamma(W,\omega^p_{D_j}) \cap L^{p,0}_\sigma(W,L_{D_j})\ ,\ \ j=1,2.$$
With $i_0=i_*|_{\ker \dq^{p,1}_{D_1}(W)}$, we see that
\begin{eqnarray}\label{eq:i0}
\Im i_0=\Im i_* \cap \ker \dq^{p,1}_{D_2}(W).
\end{eqnarray}
Then the codimension of $\Im i_0$ in $\ker\dq^{p,1}_{D_2}(W)$ is $\leq \mbox{codim}\ \Im i_*$.
This follows from the simple fact that when $U, V$ are subspaces of a vector space $W$, 
then the natural map $\frac{U}{U\cap V}\rightarrow \frac{W}{V}$, induced by the inclusion $U\hookrightarrow W$,
is injective.
Thus $\Im i_0$ is of finite codimension in $\ker\dq_{D_2}^{p,1}(W)$,
hence closed by Lemma \ref{lem:fa2}, and $i_0^{-1}: f\mapsto (f\cdot s_{D_2}^{-1})\otimes s_{D_1}$, $f\in \Im i_0$,
is bounded.

Setting $H_1'':=\{f\in H_1': \left.{\bf S}_2\circ\Psi_1(f)\right|_W \in \Im i_0\}$,
we see that $H_1''$ is closed and of finite codimension in $H_1'$.
Moreover, we see that $f\mapsto \big({\bf S}_2\circ\Psi_1(f)\cdot s^{-1}_{D_2}\big)\otimes s_{D_1}$
is a bounded linear $\dq_{D_1}$-solution operator $H_2''\rightarrow L^{p,0}_\sigma(U,L_{D_1})$
that is bounded if ${\bf S}_2$ is compact. 
%As above, $H_2''\subset \Im\dq_{D_1}$ implies that
%$\Im\dq_{D_1}$ is closed of finite codimension in $\ker\dq_{D_1}$.
\end{proof}

By a slight modification of our arguments above,
we can drop the condition on finite codimension and obtain:

\begin{thm}\label{thm:range1b}
Let $U\subset\subset M$ be a neighborhood of the exceptional set $E$ and $D_1$, $D_2$ two divisors
with support on $E$. Then
\begin{eqnarray}\label{eq:dq1b}
\dq_{D_1}: L^{p,q}_{\sigma}(U,L_{D_1}) \rightarrow L^{p,q+1}_{\sigma}(U,L_{D_1})
\end{eqnarray}
has closed range in $\ker\dq_{D_1}$ exactly if
\begin{eqnarray}\label{eq:dq2b}
\dq_{D_2}: L^{p,q}_{\sigma}(U,L_{D_2}) \rightarrow L^{p,q+1}_{\sigma}(U,L_{D_2})
\end{eqnarray}
has closed range in $\ker \dq_{D_2}$. 

If this is the case, then there exists a compact
$\dq$-solution operator 
\begin{eqnarray}\label{eq:dq11b}
{\bf S}_1: \Im \dq_{D_1} \subset L^{p,q+1}_{\sigma}(U,L_{D_1}) \rightarrow L^{p,q}_{\sigma}(U,L_{D_1})
\end{eqnarray}
exactly if there exists a compact $\dq$-solution operator
\begin{eqnarray}\label{eq:dq12b}
{\bf S}_2: \Im \dq_{D_2} \subset L^{p,q+1}_{\sigma}(U,L_{D_2}) \rightarrow L^{p,q}_{\sigma}(U,L_{D_2}).
\end{eqnarray}
\end{thm}

\begin{proof}
We just need to modify the proof of Theorem \ref{thm:range1} at some steps.

Assume first that $\dq_{D_2}$ has closed range.
We have that $H_1''$ is closed and of finite codimension in $H_1'$, and that $H_1''\subset \Im\dq_{D_1}$.
(These facts did not depend on finite codimension of $\Im \dq_{D_2}$ in $\ker\dq_{D_2}$.)

On the other hand, let $f\in\Im\dq_{D_1}$ with $f=\dq_{D_1}u$.
As $\Im\dq_{D_1} \subset H_1$,
we see (using $D_2\geq D_1$) that $\Psi_1(f)=\dq_{D_2} v$ where
$$v=\left(\big(u-\chi {\bf T}_1 (f|_V)\big)\cdot s_{D_1}^{-1}\right)\otimes s_{D_2} \in L^{p,q}_\sigma(U,L_{D_2}).$$
So, $f\in H_1'$. Hence $\Im\dq_{D_1} \subset H_1'$ and $H_1''$ has finite codimension in $\Im\dq_{D_1}$.
But then $\Im\dq_{D_1}$ is closed by Lemma \ref{lem:fa1}.

If, moreover, ${\bf S}_2: \Im\dq_{D_2} \rightarrow L^{p,q}_\sigma(U,L_{D_2})$ 
is a compact solution operator on $\Im \dq_{D_2}$, then ${\bf S}_1'$ is a compact solution operator
on $H_1''$ (see \eqref{eq:complement3}).
Hence, ${\bf S}_1: \Im \dq_{D_1} \rightarrow L^{p,q}_\sigma(U,L_{D_1})$ is compact as well by Lemma \ref{lem:fa3}
as $H_1''$ has finite codimension in $\Im\dq_{D_1}$ (see \eqref{eq:complement4} and the arguments thereafter).

\medskip

For the converse direction of the statement, assume that $\dq_{D_1}$ has closed range.
Then we know that $H_2'=\Psi_2^{-1}(\Im \dq_{D_1})$ is closed and contained in $\Im\dq_{D_2}$.
We shall find a subspace $\wt{H}_2$ of finite codimension in $\Im \dq_{D_2}$
with $\wt{H}_2 \subset H_2'$, so that $H_2'$ must have finite codimension in $\Im\dq_{D_2}$
and Lemma \ref{lem:fa1} can be applied.

Choose a (possibly discontinuous) $\dq_{D_2}$-solution operator
$${\bf S}_0: \Im\dq_{D_2} \rightarrow \Dom(\dq_{D_2}).$$
When $q>0$, consider the linear map
$$F: \Im \dq_{D_2} \rightarrow \ker\dq_{D_2}|_V / \Im \dq_{D_2}|_V \cong \C^L$$
induced by $f \mapsto ({\bf S}_0 f)|_V - {\bf T}_2 (f|_V)$.
Then $\wt{H}_2:=F^{-1}(0)$ has finite codimension in $\Im\dq_{D_2}$.
Let $f\in \wt{H}_2$. Then 
$$({\bf S}_0f)|_V-{\bf T}_2(f|_V)=\dq_{D_2} v$$ 
with $v\in L^{p,q-1}_\sigma(V,L_{D_2})$.
Observe that
\begin{eqnarray*}
f-\dq_{D_2}\big(\chi {\bf T}_2(f|_V)\big) &=& \dq_{D_2}\big((1-\chi){\bf S}_0 f\big) + \dq_{D_2}\big(\chi({\bf S}_0 f - {\bf T}_2(f|_V))\big)\\
&=& \dq_{D_2} \big((1-\chi) {\bf S}_0 f - \dq \chi\wedge v\big),
\end{eqnarray*}
since $\dq_{D_2}(\dq\chi\wedge v)= - \dq_{D_2}(\chi\dq_{D_2}v)$. 
Then $\Psi_2(f)=\dq_{D_1}w$, where
$$w=\big((1-\chi) {\bf S}_0 f -\dq\chi\wedge v\big)\cdot s_{D_2}^{-1} \otimes s_{D_1} \in L^{p,q}_\sigma(U,L_{D_1}),$$
so that $\wt{H}_2 \subset H_2'$ as desired.

\medskip

It remains to treat the case $q=0$.
For that, recall the construction of the maps $i_*$ and $i_0$ in the last section of the proof of Theorem \ref{thm:range1},
but replace the neighborhood $W$ of the exceptional set by the bigger neighborhood $V$ which has strongly pseudoconvex boundary,
as well.
We denote the corresponding maps by $i_*^V$ and $i_0^V$, respectively, i.e.
$$i_0^V = i_*^V\big|_{\ker\dq_{D_1}^{p,1}(V)}.$$

Consider now the linear map
$$G: \Im \dq_{D_2} \rightarrow \ker\dq_{D_2}(V) / \Im i_0^V \cong \C^L$$
induced by $f \mapsto ({\bf S}_0 f)|_V - {\bf T}_2 (f|_V)$.
Recall that we showed that $\ker \dq_{D_2}(V)/\Im i_0^V$ is of finite dimension.
It follows that $\wt{H}_2:=G^{-1}(0)$ has finite codimension in $\Im\dq_{D_2}$.

For $f\in \wt{H}_2$, observe as above that $\Psi_2(f)=\dq_{D_1} w$ for
$$w=\big((1-\chi){\bf S}_0 f + \chi( {\bf S}_0 f - {\bf T}_2(f|_V))\big) \cdot s^{-1}_{D_2} \otimes s_{D_1} \in L^{p,0}_\sigma(U,L_{D_1}).$$
Hence $\wt{H}_2\subset H_2'$ as desired.\\

Assume moreover that ${\bf S}_1: \Im\dq_{D_1} \rightarrow L^{p,q}_\sigma(U,L_{D_1})$ 
is a compact solution operator on $\Im \dq_{D_1}$.
Then we have seen in the proof of Theorem \ref{thm:range1} that
${\bf S}_2'$ is a compact solution operator on $H_2'$.
Hence, ${\bf S}_2: \Im \dq_{D_2} \rightarrow L^{p,q}_\sigma(U,L_{D_2})$ is compact as well by Lemma \ref{lem:fa3}
as $H_2'$ has finite codimension in $\Im\dq_{D_2}$ (see also Lemma \ref{lem:fa5}).
\end{proof}

\medskip

\section{Regularity of the $\dq$-operator at isolated singularities}\label{sec:sop}

In this section, we recall some $L^2$-regularity results for the $\dq$-equation at isolated singularities
due to Forn{\ae}ss, {\O}vrelid and Vassiliadou (see \cite{FOV2}) which imply the closed range property
of the $\dq$-operator at isolated singularities, and we deduce the existence of local compact solution operators.

As always, let $(X,g)$ be a Hermitian complex space and let $a\in X$ be an isolated singularity.
Let $U$ be a neighborhood of the point $a$ so that there exists a holomorphic embedding of $U$ in a domain $D\subset\subset \C^N$
such that $a=0$ and the metric $g$ is the pull-back of a smooth Hermitian metric $G$ on $D$.
For the following estimates, we can assume that $G$ is just the Euclidean metric, $G=i\partial\dq \|z\|^2$.
Then we can use the following results of Forn{\ae}ss, {\O}vrelid and Vassiliadou (see \cite{FOV2}) on
regularity of the $\dq$-equation at isolated singularities.

\begin{thm}\label{thm:fov2}
Let $X$ be a pure $n$-dimensional complex analytic set in $\C^N$ with an isolated singularity at $0$
and let $X$ carry the pull-back of the Euclidean metric.
Let $U$ be the intersection of $X$ with a small ball centered at the origin, $U=X\cap B_r(0)$,
and $U^*=\Reg U = U\setminus\{0\}$.

If $p+q<n$, $q\geq 1$,
then the $\dq$-operator 
$$\dq: L^{p,q-1}(U^*) \rightarrow L^{p,q}(U^*),$$
has closed image $\Im \dq$ of finite codimension in $\ker\dq$,
and there exists a constant $C>0$ such that for each $f\in \Im \dq$ there exists $u\in L^{p,q-1}(U^*)$
with $\dq u=f$ satisfying
\begin{eqnarray}\label{eq:est1}
\int_{U^*} |u|^2 \|z\|^{-2} \log^{-4} (\|z\|^{-2}) dV_X \leq C \int_{U^*} |f|^2 dV_X.
\end{eqnarray}
Let $0<c<1$.
If $p+q>n$, there exist a constant $C_c>0$ 
such that 
for each $f\in \ker\dq \subset L^{p,q}(U^*)$ there exists $u\in L^{p,q-1}(U^*)$
with $\dq u=f$ satisfying
\begin{eqnarray}\label{eq:est2}
\int_{U^*} |u|^2 \|z\|^{-2c} dV_X \leq C_c \int_{U^*} |f|^2 dV_X.
\end{eqnarray}
\end{thm}

\begin{proof}
Let $p+q<n$, and let $B$ be the $L^2$-space of $(p,q-1)$-forms $v$ such that 
$$\int_{U^*} |v|^2 \|z\|^{-2}\log^{-4} ( \|z\|^{-2}) dV_X < \infty.$$

By \cite{FOV2}, Theorem 1.1, we know that the range of $\dq: L^{p,q-1}(U^*) \rightarrow L^{p,q}(U^*)$
contains a closed subspace that has finite codimension in $\ker\dq \subset L^{p,q}(U^*)$.
Thus, Lemma \ref{lem:fa1} tells us that $\Im \dq$ is also closed and of finite codimension in $\ker\dq$.

For $f\in\Im\dq$ we choose $u\in L^{p,q-1}(U^*)$ such that $\dq u=f$.
Let $\chi \in C^\infty_{cpt}(U^*)$ be a smooth cut-off function with compact support that is identically $1$ in a neighborhood
of the origin. Then
$$f=\dq u = \dq\big((1-\chi) u\big ) +\dq\chi\wedge u + \chi f,$$
and
$$f_0:=\dq\chi\wedge u +\chi f$$ 
is $\dq$-closed with compact support in $U$.
By \cite{FOV2}, Proposition 3.1, there exists $u_0\in B$ such that $\dq u_0 = f_0$.
Thus $u':=(1-\chi)u + u_0 \in B$ and $\dq u'=f$.

That shows the existence of the desired solution $u'\in B$. The constant $C>0$ is obtained by the oppen mapping theorem:
as the mapping $\dq: \Dom\dq \subset B \rightarrow \Im \dq$ is surjective and $\Im\dq$ is closed,
there exists a continuous inverse mapping.

In the case $p+q>n$, the statement follows directly from Theorem 5.9 in \cite{FOV2}.
Here, the $\dq$-equation is solvable for all $f\in\ker\dq$.
\end{proof}

We are now in the position to construct compact solution operators for the $\dq$-equation.
Let $a\in\Sing X$ be an isolated singularity and $U$ a strongly pseudoconvex neighborhood of $a$
as in Theorem \ref{thm:fov2}.

We choose compact neighborhoods $K_1, K_2$ of the singularity $a$, with $K_1 \subset \mathring{K_2}$.
Set $U_j= U\setminus K_j$ for $j=1, 2$, and $V_1:= \mathring{K_2}\setminus K_1$.
Choose a smooth weight function $\varphi$ on $U^*$ such that $\varphi=\log \|z\|$ near $a=0$ and $\varphi\equiv 0$ on $U_1$.

Recall that $L^{p,q}(U^*,\varphi)$ is the Hilbert space of forms $f$ that are $L^2$ with respect
to the weight $\varphi$ in the sense that
$$\int_{U^*} |f|^2 e^{-\varphi} dV <\infty.$$

For $p+q\neq n$, $q\geq 1$, let
$$T: L^{p,q-1}(U^*, \varphi) \rightarrow L^{p,q}(U^*)$$
and
$$T': L^{p,q}(U^*) \rightarrow L^{p,q+1}(U^*)$$
be the $\dq$-operators in the sense of distributions.

$T$ and $T'$ are closed densely defined operators, $T' \circ T=0$ and
$T$ has closed range $\Ra (T)$ of finite codimension in $\ker T' \subset L^{p,q}(U^*)$ by Theorem \ref{thm:fov2}
(we have chosen a suitable weight). We denote the image of $T$ in $L^{p,q}(U^*)$ by $H$.
Note that $H=\ker T'$ if $p+q>n$.
So, the adjoint operators $T^*$ and $(T')^*$ are closed densely defined operators with
$T^*\circ (T')^*=0$ and $T^*$ has closed range (see \cite{Hoe1}, Theorem 1.1.1). 
Let 
\begin{eqnarray}%\label{eq:S}
{\bf S}: H \rightarrow (\ker T)^\perp= \Ra(T^*) \subset L^{p,q-1}(U^*,\varphi)
\end{eqnarray}
be the minimal solution operator for the $\dq$-equation
(see Section \ref{ssec:linops}) and extend it to an operator 
${\bf S}: L^{p,q}(U^*) \rightarrow L^{p,q-1}(U^*,\varphi)$
by setting ${\bf S}f=0$ for $f\in\ker T^*$.

\medskip

Since $L^{p,q-1}(U^*,\varphi)$ is naturally contained in $L^{p,q-1}(U^*)$,
we can show by use of the criterion for precompactness Theorem \ref{thm:precompact1}
that ${\bf S}$ is compact as an operator to the latter space.
The key idea is to exploit the fact that ${\bf S}$ is the canonical solution operator for the $\dq$-equation.
So, it can be represented as
${\bf S} = T^* N_T$,
where $N_T$ is the $\dq$-Neumann operator for the complex Laplacian $TT^*+(T')^*T'$.

\begin{thm}\label{thm:sop}
Let $p+q\neq n$, $q\geq 1$.
The $\dq$-solution operator ${\bf S}$ defined above is compact as an operator
$${\bf S}: \Dom {\bf S}=H \subset L^{p,q}(U^*) \rightarrow L^{p,q-1}(U^*).$$
\end{thm}

\begin{proof}
Let
$$\mathcal{B} = \{f\in H: \|f\|_{L^{p,q}(U^*)} <1\}.$$
We will see that $\mathcal{K}:={\bf S}(\mathcal{B})$ is relatively compact in $L^{p,q-1}(U^*)$ by Theorem \ref{thm:precompact1}.

Recall that $T$ and $T^*$ have closed range. 
Define as in Section \ref{ssec:linops} the minimal solution operator
$${\bf S}^*: L^{p,q-1}(U^*,\varphi) = \Ra(T^*)\oplus \ker T \rightarrow L^{p,q}(U^*)$$
for the operator $T^*$ and recall that ${\bf S}^*$ is in fact the Hilbert space adjoint of ${\bf S}$ (Lemma \ref{lem:fa4})
so that $\|{\bf S}^*\|=\|{\bf S}\|$ and $\|{\bf S}^* {\bf S}\|\leq \|{\bf S}\|^2$.
We define $\Box'$ by the assignment
$$\Box' \alpha := TT^* \alpha + (T')^* T'\alpha$$
when
$$\alpha\in \Dom\Box' = \{\alpha\in\Dom T'\cap\Dom T^*: T'\alpha\in\Dom (T')^*, T^*\alpha\in\Dom T\}.$$
Let $\theta=-*\partial *$ be the formal $L^2$-adjoint of $\dq$ (without any weight).
Then we have
$$(T\alpha,\beta)_{L^{p,q}(U^*)} = (\alpha,\theta \beta)_{L^{p,q-1}(U^*)} = (\alpha,e^\varphi \theta \beta)_{L^{p,q-1}(U^*,\varphi)}$$
for test-forms $\alpha$, $\beta$. Thus $T^* \beta = e^\varphi \theta \beta$
for $\beta\in\Dom T^*$, and
$$\Box' \beta = e^\varphi \theta \dq \beta + e^\varphi \dq\theta \beta + e^\varphi \dq\varphi\wedge \theta \beta$$
in $U^*$ for $\beta \in\Dom \Box'$.

Thus $\Box'\beta=\gamma$ is a determined second order elliptic system in $U^*$,
and by standard interior elliptic regularity, we have 
\begin{eqnarray}\label{eq:ellipticest1}
\|\beta\|_{W^{2,2}(V)}\ ,\ \|T^*\beta\|_{W^{1,2}(V)} \leq C_V \big(\|\beta\|_{L^2(V)} +\|\gamma\|_{L^2(V)}\big)
\end{eqnarray}
when $V\subset\subset U^*$. For elliptic and subelliptic estimates,
we drop the degree $(p,q)$ of the forms under consideration from the notation.
We write $W^{k,2}$ for the Sobolev $L^2$-spaces and just $L^2$ instead of $W^{0,2}$.

Consider $u={\bf S} f \in \mathcal{K}={\bf S}(\mathcal{B})$.
By Proposition \ref{prop:fa7}, $v:={\bf S}^* u={\bf S}^* {\bf S} f=N' f$,
where $N'$ is the $\dq$-Neumann operator for $\Box'$.
But then $u=T^* v$ and $\Box' v=f$.
Hence \eqref{eq:ellipticest1}
and the Rellich embedding theorem (i.e. the embedding $W^{1,2}(V)\hookrightarrow L^2(V)$ is compact)
show that $\mathcal{K}|_V$ is relatively compact in $L^{p,q-1}(V)$ for any $V\subset\subset U^*$.

\medskip

We must also consider the open set $U_2$ which is not relatively compact in $U^*$.
We need some subelliptic estimates at $bU$.
For that, let
$\pi: U' \rightarrow U$
be a resolution of singularities such that $U'$ appears as a strongly pseudoconvex domain in a complex manifold,
and choose a smooth metric $\sigma$ on $U'$ that equals the pull-back
of the metric from $U$ over $U_1$.

Let $\chi\in C^\infty_{cpt}(U_1\cup bU)$ be a smooth cut-off function such that $|\chi|\leq 1$ and $\chi\equiv 1$ on $U_2$.
As before, let $v={\bf S}^*{\bf S} f$ so that $\Box' v=f$.
Then $\chi v\in \Dom\Box'$ and we set
\begin{eqnarray*}
f' := \Box'(\chi v) = \chi f + \theta (\dq\chi\wedge v) - \dq (\dq\chi \lrcorner v),
\end{eqnarray*}
as $\varphi \equiv 0 $ on $U_1$ (i.e. $e^{-\varphi}\equiv 1$ on $U_1$).
But $\supp \dq \chi \subset V_1 = \mathring{K_2}\setminus K_1 \subset \subset U^*$.
So, \eqref{eq:ellipticest1} gives
\begin{eqnarray*}
\|v\|_{W^{2,2}(V_1)} \leq C_{V_1} \big( \|v\|_{L^2(V_1)} + \|f\|_{L^2(V_1)}\big) \leq C_{V_1}(1+\|{\bf S}\|^2) \|f\|_{L^2(V_1)}.
\end{eqnarray*}
This yields (use that $\chi v=v$ on $U_2$):
\begin{eqnarray}\label{eq:ff1}
\|f'\|_{L^2(U_1)} &=& \|\Box'(\chi v)\|_{L^2(U_1)} \leq \|\Box'(\chi v)\|_{L^2(V_1)} + \|\Box' v\|_{L^2(U_2)}\\
&\leq& C_{V_1}' \|f\|_{L^2(V_1)} + \|f\|_{L^2(U_2)} \leq C' \|f\|_{L^2(U_1)}\label{eq:ff2}
\end{eqnarray}
with a constant $C'>0$.

When $g$ is a form on $U_1$, we denote by $\wt{g}$ the trivial extension by $0$ of $\pi^* g$ to $U'$.
Then we see that $\wt{\chi v} \in \Dom(\Box_{U'})$ where $\Box_{U'}$ is the usual $\dq$-Laplacian on $U'$,
and $\Box_{U'}(\wt{\chi v}) = \wt{f'}$.
As $bU'$ is strongly pseudoconvex,
Kohn's subelliptic estimates yield
\begin{eqnarray*}
\|\wt{\chi v}\|_{W^{1,2}_\sigma(U')}\ ,\  \|\dq^* (\wt{\chi v})\|_{W^{1/2,2}_\sigma(U')} \leq \wt{C} \big(\|\wt{\chi v}\|_{L^2(U')} + \|\wt{f'}\|_{L^2(U')}\big).
\end{eqnarray*}
But $\sigma$ equals the pull-back of the original metric on $U_1$ and so 
$$\dq^* (\wt{\chi v}) = \wt{T^* (\chi v)} = \pi^* u$$
over $U_2$ and
\begin{eqnarray}\label{eq:u1/2}
\|u\|_{W^{1/2,2}(U_2)} \leq \wt{C} \big(\|v\|_{L^2(U_1)} + \|f'\|_{L^2(U_1)}\big) \leq C'' \|f\|_{L^2(U_1)},
\end{eqnarray}
where we have used \eqref{eq:ff1}-\eqref{eq:ff2} for the last inequality.

So, \eqref{eq:u1/2}
and the Rellich embedding theorem (i.e. the embedding $W^{1/2,2}(U_2)\hookrightarrow L^2(U_2)$ is compact)
show that $\mathcal{K}|_{U_2}$ is relatively compact in $L^{p,q-1}(U_2)$.

\medskip

Now, finally, we can apply Theorem \ref{thm:precompact1} to show that $\mathcal{K}={\bf S}(\mathcal{B})$
is in fact relatively compact in $L^{p,q-1}(U^*)$.

We cover $U^*$ by $U_2$ and by open sets $V_\mu \subset\subset U^*$.
We have already seen that condition (i) of Theorem \ref{thm:precompact1} is fulfilled
as $\mathcal{K}|_{U_2}$ and $\mathcal{K}|_{V_\mu}$ are relatively compact in $L^{p,q-1}(U_2)$ and $L^{p,q-1}(V_\mu)$, respectively,
for all $\mu$.

As $U_2$ covers $U^*$ up to the boundary $bU$, only the singularity $a=0$ is relevant in view of condition (ii).
Since ${\bf S}$ is bounded as an operator to $L^{p,q-1}(U^*,\varphi)$, there exists a constant $C_S>0$ such that
$\|u\|^2_{L^{p,q-1}(U^*,\varphi)} \leq C_S$ for all $u\in \mathcal{K}=S(\mathcal{B})$. Let $\epsilon>0$.
Choose $U_\epsilon := U\setminus B_\delta(0)$ such that
\begin{eqnarray}\label{eq:ee}
e^{-\varphi} =\|z\|^{-1} \geq 1/\epsilon
\end{eqnarray}
on $U-U_\epsilon$. Then
\begin{eqnarray*}
\epsilon^{-1} \int_{U-U_\epsilon} |u|^2 dV 
\leq  \int_{U-U_\epsilon} |u|^2 e^{-\varphi} dV
\leq \int_{U^*} |u|^2 e^{-\varphi} dV \leq C_S
\end{eqnarray*}
for all $u\in \mathcal{K}={\bf S}(\mathcal{B})$. 

Condition (ii) of Theorem \ref{thm:precompact1} follows easily
because any set $U_\epsilon$ can be covered by $U_2$ and finitely many sets $V_\mu$.
\end{proof}

We also need to treat the case $(p,q)=(0,n)$. This situation turns out to be even simpler
as the cohomology in top degree behaves nicely on non-compact complex manifolds.

\begin{thm}\label{thm:siu}
Let $X$ be a pure $n$-dimensional complex analytic set in $\C^N$ with an isolated singularity at $0$
and let $X$ carry the pull-back of the Euclidean metric.
Let $U$ be a neighborhood of $0$ in $X$ with smooth strongly pseudoconvex boundary
such that $U^*=\Reg U = U\setminus\{0\}$.
Then the $\dq$-operator 
$\dq^{0,n}: L^{0,n-1}(U^*) \rightarrow L^{0,n}(U^*)$
satisfies $\Im \dq^{0,n} = L^{0,n}(U^*)$ and there exists a compact operator
$${\bf S}: L^{0,n}(U^*) \rightarrow \Dom(\dq^{0,n}) \subset L^{0,n-1}(U^*)$$
such that $\dq^{0,n}{\bf S} f=f$ for all $f\in L^{0,n}(U^*)$.
\end{thm}

\begin{proof}
Let
$\pi: M \rightarrow X$
be a resolution of singularities as in Section \ref{sec:resolution}, $U'=\pi^{-1}(U)$ and $E=\pi^{-1}(\{0\})$.
Let $\sigma$ be any (positive definite) Hermitian metric on $M$.
We denote by $L^{p,q}$ the spaces of $L^2$-forms on open subsets of $\Reg X$,
and by $L^{p,q}_\sigma$ the spaces of $L^2$-forms on open subsets of $M$.
Moreover, let $\gamma:=\pi^* h$ be the pullback of the Hermitian metric of $X$ to $M$.
$\gamma$ is positive semidefinite (a pseudo-metric) with degeneracy locus $E$.
We denote by $L^{p,q}_\gamma$ the space of forms which are $L^2$ on $M$
with respect to the pseudo-metric $\gamma$.

For an open set $U\subset M$ and all $0\leq q\leq n$, there are bounded inclusions
\begin{eqnarray}\label{eq:l2est3}
L^{n,q}_{\gamma}(U) \subset L^{n,q}_{\sigma}(U)
\end{eqnarray}
and
\begin{eqnarray}\label{eq:l2est4}
L^{0,q}_{\sigma}(U) \subset L^{0,q}_{\gamma}(U).
\end{eqnarray}
For a proof, see e.g. \cite{Rp8}, Section 2.2, or \cite{PS1}, (1.5) and (1.6).
So, we can simply use the bounded inclusions
$$\pi|_{U'-E}^*: L^{0,n}(U^*) \overset{\cong}{\longrightarrow} L^{0,n}_\gamma(U') \cong L^{0,n}_\sigma(U')$$
(extend the forms trivially over $E$) and 
$$(\pi|_{U'-E}^{-1})^*: L^{0,n-1}_\sigma(U') \subset L^{0,n-1}_\gamma(U') \overset{\cong}{\longrightarrow} L^{0,n-1}(U^*).$$

The other tool that we need is Siu's vanishing theorem for the cohomology groups of dimension $n$
on non-compact complex manifolds of dimension $n$ (see \cite{Siu00}, 
but also \cite{M} where the statement was proved before for locally free coherent analytic sheaves):
If $N$ is an $n$-dimensional non-compact complex manifold and $\mathcal{F}$ a coherent analytic sheaf on $N$,
then $H^n(N,\mathcal{F})=0$. It follows that\footnote{As $U'$ is str. pseudoconvex,
we can also use Serre duality: $H^n(U',\OO)\cong H^0_{cpt}(U',\Omega^n)=0$.}
$H^{0,n}(U') \cong H^n(U',\OO) =0$.
On the other hand, it is well-known that the $L^2$- and the $L^{2,loc}$-Dolbeault cohomology coincide
on a domain with smooth strongly pseudoconvex boundary in a complex manifold (see e.g. \cite{LM}, Theorem VIII.4.1).
Thus
$H^{0,n}_{(2),\sigma}(U') =0$,
and strong pseudoconvexity of $U'$ implies (by Kohn's subelliptic estimates) the existence of a compact $\dq$-solution operator
$${\bf S}_M: L^{0,n}_\sigma(U') \rightarrow \Dom(\dq) \subset L^{0,n-1}_\sigma(U').$$
It follows that
${\bf S} := (\pi|_{U'-E}^{-1})^* \circ {\bf S}_M \circ \pi|_{U'-E}^*$
is the desired compact solution operator for $\dq^{0,n}$ on $U^*$.
\end{proof}

\bigskip

%\newpage
\section{Global regularity of the $\dq$-operator}\label{sec:global}

Let $X$ be a Hermitian complex space of pure dimension $n$,
i.e. a reduced complex space with a Hermitian metric which extends smoothly to the singular set, and let $\pi: M \rightarrow X$
be a resolution of singularities as in Section \ref{sec:resolution}.
Let $\sigma$ be any (positive definite) Hermitian metric on $M$.
We denote by $L^{p,q}$ the spaces of $L^2$-forms on open subsets of $\Reg X$,
and by $L^{p,q}_\sigma$ the spaces of $L^2$-forms on open subsets of $M$.

Let $\Omega\subset\subset X$ be a relatively compact open subset of $X$ such that the boundary of $\Omega$
does not intersect the singular set of $X$, $b\Omega\cap \Sing X=\emptyset$,
and that $\Omega$ contains only isolated singularities.
Let $\Omega^*:=\Omega-\Sing X$ and $\Omega':=\pi^{-1}(\Omega)$.

\begin{thm}\label{thm:compact}
Let $q\geq 1$ and either $p+q\neq n$ or $(p,q)=(0,n)$. Under the assumptions above,
the $\dq$-operator in the sense of distributions
\begin{eqnarray*}
\dq: L^{p,q-1}(\Omega^*) \rightarrow L^{p,q}(\Omega^*)
\end{eqnarray*}
has closed range (of finite codimension) in $\ker\dq\subset L^{p,q}(\Omega^*)$
exactly if the $\dq$-operator in the sense of distributions
\begin{eqnarray*}
\dq_M: L^{p,q-1}_\sigma(\Omega') \rightarrow L^{p,q}_\sigma(\Omega')
\end{eqnarray*}
has closed range (of finite codimension) in $\ker\dq_M \subset L^{p,q}_\sigma(\Omega')$.

If this is the case,
then there exists a compact $\dq$-solution operator
\begin{eqnarray*}
{\bf S}: \Im \dq \subset L^{p,q}(\Omega^*) \rightarrow L^{p,q-1}(\Omega^*)
\end{eqnarray*}
exactly if there exists a compact $\dq$-solution operator
\begin{eqnarray*}
{\bf S}_M: \Im\dq_M \subset L^{p,q}_\sigma(\Omega') \rightarrow L^{p,q-1}_\sigma(\Omega').
\end{eqnarray*}
\end{thm}

The phrase 'of finite codimension' is optional and may or may not be included in the statement of the theorem,
just as desired.

\begin{proof}
The strategy of the proof is similar to the proof of Theorem \ref{thm:range1}.

Let $E=\pi^{-1}(\Sing X)$ be the exceptional set of the resolution $\pi: M\rightarrow X$.
We may assume that $E$ is a divisor with only normal crossings,
i.e. the irreducible components of $E$ are regular and meet complex transversely.
However, this assumption is not necessary.
In the following, we can assume that $E=\pi^{-1}(\Sing X\cap \Omega)$
such that $\Omega'=\pi^{-1}(\Omega)$ is a neighborhood of the exceptional set $E$.
We denote by $\{a_1, ..., a_k\}$ the isolated singularities in $\Omega$,
so that the exceptional set consists of the components $E_\mu=\pi^{-1}(\{a_\mu\})$, $\mu=1, ..., k$,
which are pairwise disjoint. Note that $E_\mu$ will consist again of finitely many pairwise disjoint components
if $X$ is not irreducible at $a_\mu$.

Let $\gamma:=\pi^* h$ be the pullback of the Hermitian metric $h$ of $X$ to $M$.
$\gamma$ is positive semidefinite (a pseudo-metric) with degeneracy locus $E$.
We denote by $L^{p,q}_\gamma$ the space of forms which are $L^2$ on $M$
with respect to the pseudo-metric $\gamma$.
As in \eqref{eq:inclusion4}, fix a positive integer
$m$ such that the effective divisor $mE$ satisfies:
\begin{eqnarray}\label{eq:inclusion5}
L^{p,q}_{\sigma}(U,L_{-mE}) \subset L^{p,q}_{\gamma}(U) \subset L^{p,q}_\sigma(U,L_{mE})
\end{eqnarray}
for all $0\leq p,q\leq n$ and open sets $U\subset M$.

By Theorem \ref{thm:range1} and Theorem \ref{thm:range1b}, we can now replace 
the conditions on $\dq_M$ by conditions on the $\dq$-operator in the
sense of distributions for $L^2$-forms with values in the holomorphic line bundle $L_{-mE}$
which we denote by $\dq_{-mE}$:
\begin{eqnarray}\label{eq:dq21}
\dq_{M}: L^{p,q-1}_{\sigma}(\Omega') \rightarrow L^{p,q}_{\sigma}(\Omega')
\end{eqnarray}
has closed range (of finite codimension) in $\ker\dq_{M}\subset L^{p,q}_\sigma(\Omega')$ exactly if
\begin{eqnarray}\label{eq:dq22}
\dq_{-mE}: L^{p,q-1}_{\sigma}(\Omega',L_{-mE}) \rightarrow L^{p,q}_{\sigma}(\Omega',L_{-mE})
\end{eqnarray}
has closed range (of finite codimension) in $\ker \dq_{-mE}\subset L^{p,q}_\sigma(\Omega',L_{-mE})$.
If this is the case, then there exists a compact
$\dq$-solution operator 
\begin{eqnarray}\label{eq:dq31}
{\bf S}_M: \Im \dq_{M} \subset L^{p,q}_{\sigma}(\Omega') \rightarrow L^{p,q-1}_{\sigma}(\Omega')
\end{eqnarray}
exactly if there exists a compact $\dq$-solution operator
\begin{eqnarray}\label{eq:dq32}
{\bf S}_{-mE}: \Im \dq_{-mE} \subset L^{p,q}_{\sigma}(\Omega',L_{-mE}) \rightarrow L^{p,q-1}_{\sigma}(\Omega',L_{-mE}).
\end{eqnarray}

%\smallskip

Let $U_1, ..., U_k$ be strongly pseudoconvex neighborhoods of the isolated singularities
$a_1, ..., a_k$ such that Theorem \ref{thm:fov2}, Theorem \ref{thm:sop} and Theorem \ref{thm:siu} 
are valid on $U_\mu^*=U_\mu-\{a_\mu\}$, $\mu=1, ..., k$.
Let $U:=\bigcup_\mu U_\mu$ and $U^*=U-\Sing X$. 

By use of Theorem \ref{thm:fov2} and  Theorem \ref{thm:sop} or Theorem \ref{thm:siu}, respectively,
$\Im \dq|_{U^*}$ is closed and of finite codimension in $\ker\dq|_{U^*} \subset L^{p,q}(U^*)$
and there exists a compact $\dq$-solution operator
\begin{eqnarray*}
{\bf T}: \Im\dq|_{U^*} \subset L^{p,q}(U^*) \rightarrow L^{p,q-1}(U^*).
\end{eqnarray*}

As in Section \ref{sec:bundles} (see Theorem \ref{thm:range1} and Theorem \ref{thm:range1b}),
we will first prove the theorem including the statement about finite codimension,
and will show afterwards how the statement about finite codimension can be dropped.\\

Assume first that $\dq_M$ and thus also $\dq_{-mE}$ have closed range of finite codimension,
and let
\begin{eqnarray}\label{eq:dqmE}
{\bf S}_{-mE}: \Im \dq_{-mE} \subset L^{p,q}_\sigma(\Omega',L_{-mE}) \rightarrow L^{p,q-1}_\sigma(\Omega',L_{-mE})
\end{eqnarray}
be a corresponding bounded $\dq$-solution operator.

%We treat the case $p+q\neq n$. The case $(p,q)=(0,n)$ is easier and will be considered later.
As in the proof of Theorem \ref{thm:range1}, consider the bounded linear map
$$\Phi: \ker \dq\subset L^{p,q}(\Omega^*) \rightarrow \ker\dq|_{U^*}$$
given by $\Phi(f):=f|_{U^*}$.
Since $\Im\dq|_{U^*}$ is a closed subspace of finite codimension in $\ker\dq|_{U^*}$,
the same holds for
$$H:=\Phi^{-1}\big(\Im\dq|_{U^*}\big) = \big\{ f\in \ker \dq: f|_{U^*} \in \Im\dq|_{U^*}\big\}$$
in $\ker \dq$. 
Choose a smooth cut-off function $\chi$ which has compact support in $U=\bigcup_\mu U_\mu$
and is identically $1$ in a neighborhood of the isolated singularities $a_1, ..., a_k$.
Then we can define a bounded linear map
$$\Psi: H \rightarrow \ker \dq_{-mE} \subset L^{p,q}_{\sigma}(\Omega',L_{-mE})$$
by the assignment
$$f\in H \mapsto \big( \pi|_{\Omega'-E}^* \big(f - \dq ( \chi {\bf T} ( f|_{U^*}))\big) \big)\otimes s_{-mE}$$
since $f - \dq \big( \chi {\bf T} ( f|_{U^*})\big)$
is identically zero in a neighborhood the singularities $a_1, ..., a_k$.
By assumption, $\dq_{-mE}$ has closed range $\Im \dq_{-mE}$ of finite codimension in $\ker\dq_{-mE}$
so that
$$H':=\Psi^{-1} ( \Im \dq_{-mE})$$
is a closed subspace of finite codimension in $H$. 
As we have already seen that $H$ in turn is closed of finite codimension in $\ker\dq$,
it follows that $H'$ is a closed subspace of finite codimension in $\ker\dq$.

On the other hand, since $\Psi(H') \subset \Im\dq_{-mE}$,
we can define by use of \eqref{eq:dqmE} a $\dq$-solution operator
\begin{eqnarray*}
{\bf S}': H' \subset L^{p,q}(\Omega^*) \rightarrow L^{p,q-1}(\Omega^*)
\end{eqnarray*}
by setting
\begin{eqnarray*}
{\bf S}' (f):= (\pi|_{\Omega'-E}^{-1})^*\big(({\bf S}_{-mE}\circ\Psi(f))\cdot s_{-mE}^{-1}\big) + \chi {\bf T}(f|_{U^*}).
\end{eqnarray*}
Here, we use the natural injection $L^{p,q-1}_\sigma(\Omega',L_{-mE}) \subset L^{p,q-1}_\gamma(\Omega')$ 
yielding that 
$$(\pi|_{\Omega'-E}^{-1})^*: L^{p,q-1}_\sigma(\Omega',L_{-mE}) \rightarrow L^{p,q-1}(\Omega^*)$$
is a bounded linear map (see \eqref{eq:inclusion5}, \eqref{eq:inclusion4}).
Since we consider the $\dq$-operator in the sense of distributions on $\Omega^*$
and $\pi|_{\Omega'-E}$ is a biholomorphism, it follows that
$\dq {\bf S}' (f) = f$ for all $f\in H'$. Hence $H'\subset \Im\dq\subset\ker\dq$
so that $\Im\dq$ is closed and of finite codimension in $\ker\dq$.
Let ${\bf S}$ be an extension of ${\bf S}'$ to $\Im\dq$ by use of Lemma \ref{lem:fa5}.

If, in addition, ${\bf S}_M$ is compact, then ${\bf S}_{-mE}$ is compact by Theorem \ref{thm:range1} (or Theorem \ref{thm:range1b}, respectively).
This yields compactness of ${\bf S}$, because ${\bf T}$ is compact, as well.

%This settles the case $p+q\neq n$. Let now be $(p,q)=(0,n)$.
%Here, we can simply use the bounded inclusions (see \eqref{eq:l2est3} and \eqref{eq:l2est4})
%$$\pi|_{\Omega'-E}^*: L^{0,n}(\Omega^*) \rightarrow L^{0,n}_\sigma(\Omega') \cong L^{0,n}_\gamma(\Omega')$$
%(extend the forms trivially over $E$) and 
%$$(\pi|_{\Omega'-E}^{-1})^*: L^{0,n-1}_\sigma(\Omega') \subset L^{0,n-1}_\gamma(\Omega') \rightarrow L^{0,n-1}(\Omega^*).$$
%Hence,
%$${\bf S}' := (\pi|_{\Omega'-E}^{-1})^* \circ {\bf S}_M \circ \pi|_{\Omega'-E}^*$$
%is a bounded $\dq$-solution operator on $H':=(\pi|_{\Omega'-E}^*)^{-1} (\Im \dq_M)$ which is closed and of finite codimension
%in $\ker\dq$. We can complement ${\bf S}'$ to a bounded solution operator ${\bf S}$ on $\Im\dq$ 
%as above (see \eqref{eq:complement21} and \eqref{eq:complement22}),
%and it is clear that ${\bf S}'$ and ${\bf S}$ are compact if ${\bf S}_M$ is compact.\\

\medskip

For the converse direction of the statement, assume now that $\Im\dq$ has closed range of finite codimension
in $\ker\dq\subset L^{p,q}(\Omega^*)$, an let
\begin{eqnarray}\label{eq:dqS}
{\bf S}: \Im\dq\subset L^{p,q}(\Omega^*) \rightarrow L^{p,q-1}(\Omega^*)
\end{eqnarray}
be a corresponding bounded $\dq$-solution operator. As in \eqref{eq:dq21} -- \eqref{eq:dq32},
if follows from Theorem \ref{thm:range1} and Theorem \ref{thm:range1b} that it is enough to prove the statements for $\dq_{mE}$
and ${\bf S}_{mE}$ (instead for $\dq_M$ and ${\bf S}_M$).
%For this direction, we can treat the two cases $p+q\neq n$ and $(p,q)=(0,n)$ together.

Let $U_1, ..., U_k$ be the strongly pseudoconvex neighborhoods of the isolated singularities
as above and let $U_\mu':=\pi^{-1}(U_\mu)$, $\mu=1, ..., k$. Set $U':=\bigcup_\mu U'_\mu=\pi^{-1}(U)$.
As in the proof of Theorem \ref{thm:range1},
the operator
$$\dq_{mE}|_{U'}: L^{p,q-1}_\sigma(U',L_{mE}) \rightarrow L^{p,q}_\sigma(U',L_{mE})$$
has closed range of finite codimension in $\ker\dq_{mE}|_{U'}$ and there exists a corresponding compact $\dq$-solution operator
\begin{eqnarray}\label{eq:Tmu}
{\bf T}_{mE}: \Im\dq_{mE}|_{U'} \subset L^{p,q}_\sigma(U',L_{mE}) \rightarrow L^{p,q-1}_\sigma(U',L_{mE}).
\end{eqnarray}
Here now, let
$$\Phi_{mE}: \ker \dq_{mE} \subset L^{p,q}_\sigma(\Omega',L_{mE}) \rightarrow \ker\dq_{mE}|_{U'}$$
be the bounded linear map given by $\Phi_{mE}(f) := f|_{U'}$,
such that
$$H_{mE}:=\Phi^{-1}_{mE}(\Im\dq_{mE}|_{U'})$$
is a closed subspace of finite codimension in $\ker\dq_{mE}$.

Setting $\chi'=\pi^*\chi$, we define here a bounded linear map
$$\Psi_{mE}: H_{mE} \rightarrow \ker \dq \subset L^{p,q}(\Omega^*)$$
by the assignment
$$f\in H_{mE} \mapsto (\pi|_{\Omega'-E}^{-1})^* \big(\big(f - \dq_{mE} \chi' {\bf T}_{mE} ( f|_{U'})\big) \cdot s_{mE}^{-1}\big)$$
since $f - \dq_{mE} \chi' {\bf T}_{mE} ( f|_{U'})$
is identically zero in a neighborhood the exceptional set $E$.
By assumption, $\dq$ has closed range $\Im \dq$ of finite codimension in $\ker\dq$
so that
$$H'_{mE}:=\Psi^{-1}_{mE} ( \Im \dq)$$
has closed range of finite codimension in $H_{mE}$. 
As we have already seen that $H_{mE}$ in turn is closed of finite codimension in $\ker\dq_{mE}$,
it follows that $H'_{mE}$ is a closed subspace of finite codimension in $\ker\dq_{mE}$.

On the other hand, since $\Psi_{mE}(H'_{mE}) \subset \Im\dq$,
we can define by use of \eqref{eq:dqS} a $\dq$-solution operator
\begin{eqnarray*}
{\bf S}'_{mE}: H'_{mE} \subset L^{p,q}_\sigma(\Omega',L_{mE}) \rightarrow L^{p,q-1}_\sigma(\Omega',L_{mE})
\end{eqnarray*}
by setting
\begin{eqnarray*}
{\bf S}'_{mE} (f):= \big( \pi|_{\Omega'-E}^*({\bf S}\circ\Psi_{mE}(f)) \big)\otimes s_{mE} + \chi' {\bf T}_{mE}(f|_{U'}).
\end{eqnarray*}
Here, we use the natural injection $L^{p,q-1}_\gamma(\Omega')\subset L^{p,q-1}_\sigma(\Omega',L_{mE})$ 
yielding that 
$$\pi|_{\Omega'-E}^*: L^{p,q-1}(\Omega^*) \rightarrow L^{p,q-1}_\sigma(\Omega',L_{mE})$$
is a bounded linear map (see \eqref{eq:inclusion5}, \eqref{eq:inclusion4}).

In this situation, we only know that
$\dq_{mE} {\bf S}_{mE}' (f) = f$
on $\Omega'\setminus E$. But, the $\dq_{mE}$-equation in $L^{p,*}_\sigma(\Omega,L_{mE})$ extends over the hypersurface $E$
so that $\dq_{mE}{\bf S}_{mE}' (f)=f$ on $\Omega'$ (see the $\dq$-extension Theorem 3.2 in \cite{Rp1}).
Hence, here $H'_{mE}\subset \Im\dq_{mE}\subset\ker\dq_{mE}$
so that $\Im\dq_{mE}$ is closed and of finite codimension in $\ker\dq_{mE}$.

We can complement ${\bf S}_{mE}'$ by Lemma \ref{lem:fa5} to a bounded $\dq$-solution
operator
$${\bf S}_{mE}: \Im\dq_{mE} \subset L^{p,q}_\sigma(\Omega',L_{mE}) \rightarrow L^{p,q-1}_\sigma(\Omega',L_{mE}).$$
If in addition ${\bf S}$ is compact, then compactness of ${\bf T}_{mE}$
implies that ${\bf S}_{mE}'$ and ${\bf S}_{mE}$ are also compact.

\medskip

Finally, we will now prove that 
$\dq^{p,q}: L^{p,q-1}(\Omega^*)\rightarrow L^{p,q}(\Omega^*)$
has closed range exactly if
$\dq_M^{p,q}: L^{p,q-1}_\sigma(\Omega') \rightarrow L^{p,q}_\sigma(\Omega')$
has closed range, i.e that we can drop the statement about finite codimension.

\medskip

Observe that for $q=1$, we have injective maps
$$i_1: \ker\dq^{p,1}_{-mE}|_{U'} \rightarrow \ker\dq^{p,1}|_{U^*}\ ,\ h\mapsto (\pi^{-1}_{U'-E})^* \big(h\cdot s_{-mE}^{-1}\big)$$
and
$$i_2: \ker\dq^{p,1}|_{U^*} \rightarrow \ker\dq^{p,1}_{mE}|_{U'}\ ,\ f\mapsto (\pi^* f) \otimes s_{mE}$$
such that $i_2\circ i_1=i_0$ where
$$i_0: \ker\dq^{p,1}_{-mE}|_{U'} \rightarrow \ker\dq^{p,1}_{mE}|_{U'}\ ,\ g\mapsto (g\cdot s_{-mE}^{-1})\otimes s_{mE}.$$
Note that we have considered this map $i_0$ already in the last section of the proof of Theorem \ref{thm:range1}
and that $\Im i_0$  has finite codimension in $\ker\dq^{p,1}_{mE}|_{U'}$ (see \eqref{eq:i0} and use $D_1=-mE$, $D_2=mE$, $W=U'$).
But $\Im i_0 \subset \Im i_2$ so that $\Im i_2$ has finite codimension in $\ker\dq^{p,1}_{mE}|_{U'}$, as well.
On the other hand, $\Im i_1= i^{-1}_2(\Im i_0)$ as $i_2$ is injective.
Hence, $\Im i_1$ is of finite codimension in $\ker\dq^{p,1}|_{U^*}$.
We conclude by use of Lemma \ref{lem:fa2} that both, $\Im i_1$ and $\Im i_2$, are closed.

\medskip

Assume now that $\dq^{p,q}_{-mE}$ has closed range (which is equivalent to $\dq_M^{p,q}$ having closed range by Theorem \ref{thm:range1b}).
Then
$$H'=\Psi^{-1}(\Im \dq^{p,q}_{-mE})$$
is a closed subspace of $\Im \dq^{p,q}$. We shall find a another subspace $\wt{H}\subset H'$
of finite codimension in $\Im \dq^{p,q}$. Then we can conclude that $\Im\dq^{p,q}$
is closed by use of Lemma \ref{lem:fa1}.
As in the proof of Theorem \ref{thm:range1b},
we choose a (possibly discontinuous) solution operator
$${\bf S}_0: \Im \dq^{p,q} \rightarrow\Dom(\dq^{p,q})$$
with $\dq^{p,q} {\bf S}_0 f=f$.
When $q>1$, set
$$\wt{H}:=\{ f\in\Im \dq^{p,q}: ({\bf S}_0 f)|_{U^*} - {\bf T}(f|_{U^*}) \in \Im \dq^{p,q-1}|_{U^*}\}.$$
Thus, $\wt{H}$ is of finite codimension in $\Im\dq^{p,q}$ because $f\mapsto ({\bf S}_0 f)|_{U^*} - {\bf T}(f|_{U^*})$
is a linear map $\Im \dq^{p,q} \rightarrow \ker \dq^{p,q}|_{U^*}$
and $\Im\dq^{p,q-1}|_{U^*}$ is of finite codimension in $\ker\dq^{p,q}|_{U^*}$.
Let $f\in\wt{H}$ with
$$({\bf S}_0 f)|_{U^*} - {\bf T}(f|_{U^*}) = \dq^{p,q-1}|_{U^*} v$$
for $v\in L^{p,q-2}(U^*)$. Observe that
\begin{eqnarray*}
f-\dq^{p,q} \big(\chi {\bf T}(f|_{U^*})\big) &=& \dq^{p,q} \big((1-\chi){\bf S}_0 f\big) + \dq^{p,q}\big(\chi({\bf S}_0 f - {\bf T}(f|_{U^*}))\big)\\
&=& \dq^{p,q} \big((1-\chi) {\bf S}_0 f - \dq \chi\wedge v\big),
\end{eqnarray*}
since $\dq^{p,q}(\dq\chi\wedge v)= - \dq^{p,q}(\chi\dq^{p,q}v)$. 
It follows that $\Psi(f)=\dq^{p,q}_{-mE} w$ where
$$w= \pi^*\big( (1-\chi){\bf S}_0 f - \dq\chi\wedge v)\big) \otimes s_{-mE} \in L^{p,q-1}_\sigma(\Omega',L_{-mE})$$
so that in fact $f\in H'$. So, $\wt{H}\subset H'$ as desired.

For the case $q=1$, we set instead
$$\wt{H}=\{f\in \Im\dq^{p,1}: ({\bf S}_0 f)|_{U^*} - {\bf T}(f|_{U^*}) \in \Im i_1\}.$$
Thus, $\wt{H}$ is of finite codimension in $\Im\dq^{p,1}$ because $f\mapsto ({\bf S}_0 f)|_{U^*} - {\bf T}(f|_{U^*})$
is a linear map $\Im \dq^{p,1} \rightarrow \ker \dq^{p,1}|_{U^*}$
and $\Im i_1$ is of finite codimension in $\ker\dq^{p,1}|_{U^*}$.
Let $f\in\wt{H}$ with
$$({\bf S}_0 f)|_{U^*} - {\bf T}(f|_{U^*})=i_1(v)$$
for $v\in \ker \dq^{p,1}_{-mE}|_{U'}$. 

Observe as above that
\begin{eqnarray*}
f-\dq^{p,q} \big(\chi {\bf T}(f|_{U^*})\big) &=& \dq^{p,q} \big((1-\chi){\bf S}_0 f + \chi(({\bf S}_0 f)|_{U^*} - {\bf T}(f|_{U^*}))\big)\\
&=& \dq^{p,q} \big((1-\chi){\bf S}_0 f + \chi i_1(v)\big).
\end{eqnarray*}
It follows that $\Psi(f)=\dq^{p,q}_{-mE} u$ for
$$u= \pi^*\big( (1-\chi){\bf S}_0 f + \chi i_1 v \big) \otimes s_{-mE} \in L^{p,0}_\sigma(\Omega',L_{-mE})$$
so that in fact $f\in H'$. So, we have $\wt{H}\subset H'$ also in the case $q=1$.

\medskip

Assume moreover that ${\bf S}_{-mE}: \Im\dq_{-mE} \subset L^{p,q}_\sigma(\Omega',L_{-mE}) \rightarrow L^{p,q-1}_\sigma(\Omega',L_{-mE})$ 
is a compact solution operator on $\Im \dq_{-mE}$.
Then we have seen above that
${\bf S}': H'\subset L^{p,q}(\Omega^*) \rightarrow L^{p,q-1}(\Omega^*)$ is a compact solution operator on $H'$.
Hence, ${\bf S}: \Im \dq^{p,q} \rightarrow L^{p,q-1}(\Omega^*)$ is compact as well by Lemma \ref{lem:fa3}
as $H'$ has finite codimension in $\Im\dq^{p,q}$.

\medskip

For the converse direction of the statement, assume now that 
$$\dq^{p,q}: L^{p,q-1}(\Omega^*) \rightarrow L^{p,q}(\Omega^*)$$
has closed range $\Im\dq^{p,q}$ in $L^{p,q}(\Omega^*)$. Then we have that
$$H_{mE}'=\Psi_{mE}^{-1}(\Im\dq^{p,q})$$
is a closed subspace of $\Im \dq_{mE}^{p,q}$, and we shall find another subspace
$\wt{H}_{mE} \subset H'_{mE}$ such that $\wt{H}_{mE}$ has finite codimension in $\Im \dq_{mE}^{p,q}$.
Then $H'_{mE}$ is also of finite codimension in $\Im \dq_{mE}^{p,q}$, thus $\Im \dq_{mE}^{p,q}$ is closed
by Lemma \ref{lem:fa1}.

For that, we repeat the procedure from above.
Choose a (possibly discontinuous) solution operator
$${\bf S}_0: \Im \dq_{mE} \rightarrow \Dom(\dq^{p,q}_{mE})$$
with $\dq^{p,q}_{mE}{\bf S}_0 f=f$.
When $q>1$, set
$$\wt{H}_{mE}:=\{ f\in\Im\dq^{p,q}_{mE}: ({\bf S}_0 f)|_{U'} - {\bf T}_{mE}(f|_{U'}) \in \Im\dq^{p,q}_{mE}|_{U'}\},$$
where ${\bf T}_{mE}$ is the compact solution operator from \eqref{eq:Tmu}.
So, $\wt{H}_{mE}$ is of finite codimension in $\Im\dq^{p,q}_{mE}$.
Let $f\in\wt{H}_{mE}$ with
$$({\bf S}_0 f)|_{U'} - {\bf T}_{mE}(f|_{U'}) = \dq^{p,q-1}|_{U'} v$$
for $v\in L^{p,q-2}_\sigma(U',L_{mE})$. Observe as above that
\begin{eqnarray*}
f-\dq^{p,q}_{mE} \big(\chi' {\bf T}_{mE}(f|_{U'})\big) = \dq^{p,q}_{mE} \big((1-\chi') {\bf S}_0 f - \dq \chi'\wedge v\big).
\end{eqnarray*}
It follows that $\Psi_{mE}(f)=\dq^{p,q} w$ where
$$w= (\pi|_{\Omega'-E}^{-1})^* \left(\big( (1-\chi'){\bf S}_0 f - \dq\chi'\wedge v)\big) \cdot s_{mE}^{-1}\right) \in L^{p,q-1}(\Omega^*)$$
so that in fact $f\in H'_{mE}$. Hence, $\wt{H}_{mE}\subset H'_{mE}$ as desired.

For the case $q=1$, we set
$$\wt{H}_{mE}=\{f\in \Im\dq^{p,1}_{mE}: ({\bf S}_0 f)|_{U'} - {\bf T}(f|_{U'}) \in \Im i_2\}.$$
Thus, $\wt{H}_{mE}$ is of finite codimension in $\Im\dq^{p,1}_{mE}$ because
$\Im i_2$ is of finite codimension in $\ker\dq^{p,1}_{mE}|_{U'}$.
Let $f\in\wt{H}_{mE}$ with
$$({\bf S}_0 f)|_{U'} - {\bf T}(f|_{U'})=i_2(v)$$
for $v\in \ker \dq^{p,1}|_{U^*}$. Observe as above that
\begin{eqnarray*}
f-\dq^{p,q}_{mE} \big(\chi' {\bf T}(f|_{U'})\big) = \dq^{p,q}_{mE} \big((1-\chi'){\bf S}_0 f + \chi' i_2(v)\big).
\end{eqnarray*}
It follows that $\Psi_{mE}(f)=\dq^{p,q} u$ for
$$u= (\pi|_{\Omega'-E}^{-1})^* \left(\big( (1-\chi'){\bf S}_0 f + \chi' i_2 (v) \big)\cdot s^{-1}_{mE}\right) \in L^{p,0}(\Omega^*)$$
so that in fact $f\in H'_{mE}$. So, we have $\wt{H}_{mE}\subset H'_{mE}$ also in the case $q=1$.

\medskip

Assume moreover that ${\bf S}: \Im\dq^{p,q} \subset L^{p,q}(\Omega^*) \rightarrow L^{p,q-1}(\Omega^*)$ 
is a compact solution operator on $\Im \dq^{p,q}$.
Then we have seen above that
$${\bf S}'_{mE}: H'_{mE}\subset L^{p,q}_\sigma(\Omega',L_{mE}) \rightarrow L^{p,q-1}_\sigma(\Omega',L_{mE})$$ 
is a compact solution operator on $H'_{mE}$.
Hence, ${\bf S}_{mE}: \Im \dq^{p,q}_{mE} \rightarrow L^{p,q-1}_\sigma(\Omega',L_{mE})$ is compact as well by Lemma \ref{lem:fa3}
as $H'_{mE}$ has finite codimension in $\Im\dq^{p,q}_{mE}$.
\end{proof}

\bigskip

\section{Compactness of the $\dq$-Neumann operator}

We need to add another criterion to Theorem \ref{thm:compact3},
namely the characterization of compactness of the $\dq$-Neumann operator
by the existence of compact $\dq$-solution operators,
a criterion which holds on arbitrary Hermitian manifolds (Theorem \ref{thm:N100}).
The proof of Theorem \ref{thm:N100} is an easy consequence of the preliminaries
on closed, densely defined linear operators that we have collected in Section \ref{ssec:linops}.

Let $M$ be a Hermitian complex manifold of dimension $n$, and let $0\leq p,q\leq n$ and $q\geq 1$.
Let $\dq_{p,q}: L^{p,q-1}(M) \rightarrow L^{p,q}(M)$
be the $\dq$-operator in the sense of distributions and $\dq_{p,q}^*$ its $L^2$-adjoint.

We assume that $\dq_{p,q}$ (and so also $\dq_{p,q}^*$) have closed range and define
the minimal (i.e. canonical) solution operators for $\dq_{p,q}$ and $\dq_{p,q}^*$,
\begin{eqnarray*}
&{\bf S}_{p,q}:& \Im\dq_{p,q} \rightarrow \big(\ker \dq_{p,q} \big)^\perp = \overline{\Im\dq_{p,q}^*},\\ 
&{\bf S}^*_{p,q}:& \Im\dq_{p,q}^* \rightarrow \big(\ker\dq_{p,q}^*\big)^\perp= \overline{\Im\dq_{p,q}},
\end{eqnarray*} 
as in Section \ref{ssec:linops}. Recall that ${\bf S}_{p,q}$ and ${\bf S}_{p,q}^*$ are bounded,
and that ${\bf S}_{p,q}^*$ is in fact the $L^2$-adjoint of ${\bf S}_{p,q}$ (Lemma \ref{lem:fa4}).
Note also the trivial:

\begin{lem}\label{lem:Scpt}
${\bf S}_{p,q}$ is compact exactly if there exists a compact $\dq$-solution operator
$${\bf T}_{p,q}: \Im\dq_{p,q} \rightarrow L^{p,q-1}(M).$$
\end{lem}

\begin{proof}
Simply compose ${\bf T}_{p,q}$ with the (bounded) orthogonal projection onto $(\ker\dq_{p,q})^\perp$
and extend this operator by zero to $(\Im\dq_{p,q})^\perp$. The other direction is trivial.
\end{proof}

Now, we draw our attention to the $\dq$-Neumann operator which we will represent by use of
the canonical $\dq$-solution operators discussed above.
On
$$\Dom\Box_{p,q} = \{ u\in \Dom\dq_{p,q+1}\cap\Dom\dq^*_{p,q}: \dq_{p,q+1} f\in\Dom\dq^*_{p,q+1}, \dq^*_{p,q}f\in\Dom\dq_{p,q}\},$$
we define the $\dq$-Laplacian
$$\Box_{p,q}=\dq_{p,q}\dq^*_{p,q}+\dq^*_{p,q+1}\dq_{p,q+1}: \Dom\Box_{p,q} \subset L^{p,q}(M) \rightarrow L^{p,q}(M).$$
It is well-known that this is a densely defined, closed, self-adjoint operator (see Theorem \ref{thm:fa5}).
The minimal solution operator for $\Box_{p,q}$,
$$N_{p,q}=\Box^{-1}_{p,q}: \Im\Box_{p,q} \subset L^{p,q}(M) \rightarrow \Dom\Box_{p,q}\cap(\ker\Box_{p,q})^\perp \subset L^{p,q}(M),$$ 
is called the $\dq$-Neumann operator.
$N_{p,q}$ is a bounded operator exactly if $\Box_{p,q}$ has closed range 
which in turn is the case exactly if $\dq_{p,q}$ and $\dq_{p,q+1}$ both have closed range (see Theorem \ref{thm:fa6}).
If this is the case, we extend $N_{p,q}$ as above trivially to an operator
$$N_{p,q}: L^{p,q}(M) \rightarrow \Dom(\Box_{p,q})\cap (\ker \Box_{p,q})^\perp.$$

\begin{thm}\label{thm:N100}
Let $\dq_{p,q}$ and $\dq_{p,q+1}$ have closed range. Then:
\begin{eqnarray}\label{eq:N100}
N_{p,q} = {\bf S}_{p,q}^* {\bf S}_{p,q} + {\bf S}_{p,q+1} {\bf S}_{p,q+1}^*.
\end{eqnarray}
Hence, the $\dq$-Neumann operator $N_{p,q}$ is compact exactly if the canonical $\dq$-solution
operators ${\bf S}_{p,q}$ and ${\bf S}_{p,q+1}$ both are compact.
\end{thm}

\begin{proof}
The identity \eqref{eq:N100} is direct consequence of Theorem \ref{thm:fa6}.

For a bounded operator $T$,
it is well known that $T$ is compact exactly if $T^*$ is compact,
and this is the case exactly if $T^*T$ is compact.

So, assume that ${\bf S}_{p,q}$ and ${\bf S}_{p,q+1}$ are compact.
Then it follows from \eqref{eq:N100} that $N_{p,q}$ is compact.
Conversely, assume that $N_{p,q}$ is compact. Then ${\bf S}_{p,q}^*{\bf S}_{p,q}$
and ${\bf S}_{p,q+1} {\bf S}_{p,q+1}^*$ both are compact by \eqref{eq:N100}
for they are positive. It follows that ${\bf S}_{p,q}$ and ${\bf S}_{p,q+1}$ both
are compact.
\end{proof}

This gives another criterion for compactness of the $\dq$-Neumann operator on
arbitrary Hermitian manifolds which we apply in the context of isolated singularities
to finally prove our second main theorem.

\subsection{Proof of Theorem \ref{thm:1.2}}

Statement i. follows directly from the discussion above (use Theorem \ref{thm:fa6}).
Moreover, Theorem \ref{thm:compact} shows that
\begin{eqnarray*}
&\dq_{p,q}:& L^{p,q-1}(\Omega^*)\rightarrow L^{p,q}(\Omega^*),\\
&\dq_{p,q+1}:& L^{p,q}(\Omega^*) \rightarrow L^{p,q+1}(\Omega^*)
\end{eqnarray*}
also have closed range in the corresponding kernels of $\dq$.
So, statement ii. is also a direct consequence of the discussion above.

It follows from the proof of Theorem \ref{thm:compact} that $\Im \dq_{p,q}$ is of finite codimension in $\ker \dq_{p,q+1}$
if $\Im \dq_{p,q}^M$ is of finite codimension in $\ker \dq_{p,q+1}^M$ and that 
$\Im \dq_{p,q+1}$ is of finite codimension in $\ker \dq_{p,q+2}$
if $\Im \dq_{p,q+1}^M$ is of finite codimension in $\ker \dq_{p,q+2}^M$.
As $\Box^M_{p,q}$ and $\Box_{p,q}$ both have closed range, we have orthogonal decompositions
$L^{p,q}_\sigma(\Omega') = \Im \Box^M_{p,q} \oplus \ker \Box^M_{p,q}$ and $L^{p,q}(\Omega^*)=\Im\Box_{p,q}\oplus \ker\Box_{p,q}$.
But $\ker\Box^M_{p,q}$ and $\ker\Box_{p,q}$ are of finite dimension if $\Im \dq_{p,q}^M$ and $\Im \dq_{p,q}$
are of finite codimension (as the $L^2$-cohomology classes have harmonic representatives). That shows iii.

It remains to show that $N_{p,q}$ is compact exactly if $N_{p,q}^M$ is compact.
By Theorem \ref{thm:N100}, $N^M_{p,q}$ is compact exactly if
the canonical $\dq^M$-solution operators
\begin{eqnarray*}
{\bf S}_{p,q}^M : && L^{p,q}_\sigma(\Omega')\rightarrow L^{p,q-1}_\sigma(\Omega'),\\
{\bf S}_{p,q+1}^M: && L^{p,q+1}_\sigma(\Omega') \rightarrow L^{p,q}_\sigma(\Omega')
\end{eqnarray*}
both are compact. By Theorem \ref{thm:compact} and Lemma \ref{lem:Scpt}, 
this is exactly the case if the canonical $\dq$-solution operators
\begin{eqnarray*}
{\bf S}_{p,q} : && L^{p,q}(\Omega^*)\rightarrow L^{p,q-1}(\Omega^*),\\
{\bf S}_{p,q+1}: && L^{p,q+1}(\Omega^*) \rightarrow L^{p,q}(\Omega^*)
\end{eqnarray*}
both are compact. Another application of Theorem \ref{thm:N100}
shows that this in turn is equivalent to compactness of $N_{p,q}$.

%\bigskip


\newpage



\begin{thebibliography}{99999}

\bibitem[A]{Alt} {\sc H.\ W.\ Alt}, {\em Lineare Funktionalanalysis}, Springer-Verlag, Berlin, 1992.

\bibitem[AS1]{AS1} {\sc M.\ Andersson, H.\ Samuelsson},
A Dolbeault-Grothendieck lemma on complex spaces via Koppelman formulas,
{\em Invent. Math.} {\bf 190} (2012), 261--297. 

\bibitem[AS2]{AS2} {\sc M.\ Andersson, H.\ Samuelsson},
Weighted Koppelman formulas and the $\dq$-equation on an analytic space,
{\em J. Funct. Anal.} {\bf 261} (2011), no. 3, 777--802.

\bibitem[AV]{AnVe} {\sc A.\ Andreotti, E.\ Vesentini},
Carleman estimates for the Laplace Beltrami equation on complex manifolds,
{\em Publ. Math. Inst. Hautes Etudes Sci.} {\bf 25} (1965), 81--130.

\bibitem[AHL]{AHL} {\sc J.\ M.\ Aroca, H.\ Hironaka, J.\ L.\ Vicente},
Desingularization theorems, {\em Mem. Math. Inst. Jorge Juan}, {\bf No. 30}, Madrid, 1977.

\bibitem[BM]{BiMi} {\sc E.\ Bierstone, P.\ Milman}, Canonical desingularization in characteristic zero by blowing-up
the maximum strata of a local invariant, {\em Inventiones Math.} {\bf 128} (1997), {\em no. 2}, 207--302.

%\bibitem[CM]{CM} {\sc M.\ Coltoiu, N.\ Mihalache},
%Strongly plurisubharmonic exhaustion functions on 1-convex spaces, {\em Math. Ann.} {\bf 270} (1985), 63--68.

\bibitem[DFV]{DFV} {\sc K.\ Diederich, J.\ E.\ Forn{\ae}ss, S.\ Vassiliadou},
Local $L^2$ results for $\dq$ on a singular surface, {\em Math. Scand.} {\bf 92} (2003), 269--294.

\bibitem[FK]{FK}{\sc G.\ B.\ Folland, J.\ J.\ Kohn}, {\em The Neumann problem for the Cauchy-Riemann complex},
Ann. of Math. Stud., No. 75, Princeton University Press, Princeton, 1972.

\bibitem[F]{Fo} {\sc J.\ E.\ Forn{\ae}ss}, $L^2$ results for $\dq$ in a conic,
in {\em International Symposium, Complex Analysis and Related Topics, Cuernavaca, Operator Theory:
Advances and Applications} (Birkhauser, 1999).

\bibitem[FOV1]{FOV1} {\sc J.\ E.\ Forn{\ae}ss, N.\ {\O}vrelid, S.\ Vassiliadou},
Semiglobal results for $\dq$ on a complex space with arbitrary singularities,
{\em Proc. Am. Math. Soc.} {\bf 133} (2005), {\em no. 8}, 2377--2386.

\bibitem[FOV2]{FOV2} {\sc J.\ E.\ Forn{\ae}ss, N.\ {\O}vrelid, S.\ Vassiliadou},
Local $L^2$ results for $\dq$: the isolated singularities case,
{\em Internat. J. Math.} {\bf 16} (2005), {\em no. 4}, 387--418.

%\bibitem[FH]{FoHa}{\sc J.\ Fox, P.\ Haskell},
%Hodge decompositions and Dolbeault complexes on complex surfaces,
%{\em Trans. Amer. Math. Soc.} {\bf 343} (1994), 765--778.

%\bibitem[G1]{Gf}{\sc M.\ P.\ Gaffney},
%Hilbert space methods in the theory of harmonic integrals,
%{\em Trans. Amer. Math. Soc.} {\bf 78} (1955), 426--444.

%\bibitem[G2]{G} {\sc K.\ Gansberger},
%On the resolvent of the Dirac operator in $\R^2$,{\sf arXiv:1003.5124}.

%\bibitem[H1]{Hk}{\sc P.\ Haskell}, $L^2$-Dolbeault complexes on singular curves and surfaces,
%{\em Proc. Amer. Math. Soc.} {\bf 107} (1989), no. 2, 517--526.

%\bibitem[H2]{H} {\sc F.\ Haslinger},
%Compactness for the $\dq$-Neumann problem - a Functional Analysis Approach,
%ESI-Preprint {\bf 2208}, {\sf arXiv:0912.4406}.

\bibitem[H1]{Ha} {\sc H.\ Hauser}, The Hironaka theorem on resolution of singularities,
{\em Bull. (New Series) Amer. Math. Soc.} {\bf 40}, no. 3, (2003), 323--403.

%\bibitem[HL]{HL} {\sc T.\ Hefer, I.\ Lieb},
%On the compactness of the $\overline\partial$-Neumann operator,  
%{\em Ann. Fac. Sci. Toulouse Math. (6)}  {\bf 9}  (2000),  no. 3, 415--432. 



\bibitem[HL1]{HL} {\sc T.\ Hefer, I.\ Lieb},
On the compactness of the $\overline\partial$-Neumann operator,  
{\em Ann. Fac. Sci. Toulouse Math. (6)}  {\bf 9}  (2000),  no. 3, 415--432. 


\bibitem[HL2]{HeLe} {\sc G.\ M.\ Henkin, J.\ Leiterer}, {\em Theory of Functions on Complex Manifolds},
Birkh\"auser, Basel, 1984.




\bibitem[H2]{Hi} {\sc H.\ Hironaka}, Resolution of singularities of an algebraic variety over a field of characteristic zero: I, II. 
{\em Ann. of Math. (2)} {\bf 79} (1964), 109--326.

\bibitem[H3]{Hoe1}{\sc L.\ H\"ormander}, $L^2$-estimates and existence theorems for the $\dq$-operator,
{\em Acta Math.} {\bf 113} (1965), 89--152.

\bibitem[H4]{Hoe2}{\sc L.\ H\"ormander}, {\em An introduction to complex analysis in several variables},
Van Nostrand, Princeton, 1966.

\bibitem[H5]{Hoe3}{\sc L.\ H\"ormander}, The null space of the $\dq$-Neumann operator,
{\em Ann. Inst. Fourier (Grenoble)} {\bf 54} (2004), 1305--1369.


%\bibitem[KK]{KK}{\sc M.\ Kashiwara, T.\ Kawai}, The Poincar\'e lemma for varieties of polarized Hodge structure,
%{\em Publ. Res. Inst. Math. Sci.} {\bf 23} (1987), 345--407.

\bibitem[K1]{Ko1}{\sc J.\ J.\ Kohn}, Harmonic integrals on strongly pseudoconvex manifolds I,
{\em Ann. of Math.} {\bf 78} (1963), 112--148.

\bibitem[K2]{Ko2}{\sc J.\ J.\ Kohn}, Harmonic integrals on strongly pseudoconvex manifolds II,
{\em Ann. of Math.} {\bf 79} (1964), 450--472.

\bibitem[KN]{KoNi}{\sc J.\ J.\ Kohn, L.\ Nirenberg},
Non-coercive boundary value problems, {\em Comm. Pure Appl. Math.} {\bf 18} (1965), 443--492.

\bibitem[LM]{LM} {\sc I.\ Lieb, J.\ Michel}, 
{\em The Cauchy-Riemann Complex, Integral Formulae and Neumann Problem}, 
Vieweg, Braunschweig/Wiesbaden, 2002.



\bibitem[M]{M}{\sc B.\ Malgrange},
Faisceaux sur des vari{\'e}t{\'e}s analytiques-r{\'e}elles, {\em Bull. Soc. Math. France} {\bf 85} (1957), 231--237.



%\bibitem[N]{N2}{\sc M.\ Nagase},
%Remarks on the $L^2$-Dolbeault cohomology groups of singular algebraic surfaces and curves,
%{\em Publ. Res. Inst. Math. Sci.} {\bf 26} (1990), no. 5, 867--883.

%\bibitem[O]{Oh2}{\sc T.\ Ohsawa}, Hodge spectral sequence on compact K\"ahler spaces,
%{\em Publ. Res. Inst. Math. Sci.} {\bf 23} (1987), 265--274.

\bibitem[OT]{OT} {\sc T.\ Ohsawa, K.\ Takegoshi},
On the extension of $L^2$-holomorphic functions,
{\em Math. Z.} {\bf 195} (1987), no. 2, 197--204.

\bibitem[OV1]{OV1} {\sc N.\ {\O}vrelid, S.\ Vassiliadou}, Solving $\dq$ on product singularities,
{\em Complex Var. Ellipitic Equ.} {\bf 51} (2006), {\em no. 3}, 225--237.

\bibitem[OV2]{OV2}{\sc N.\ {\O}vrelid, S.\ Vassiliadou}, Some $L^2$ results for $\dq$ on projective varieties
with general singularities, {\em Amer. J. Math.} {\bf 131} (2009), 129--151.

\bibitem[OV3]{OV3}{\sc N.\ {\O}vrelid, S.\ Vassiliadou}, $L^2$-$\dq$-cohomology groups of some singular complex spaces,
{\em Invent. Math.} 2012, online first, DOI: 10.1007/s00222-012-0414-3.

\bibitem[P]{P}{\sc W. Pardon}, The $L^2$-$\dq$-cohomology of an algebraic surface,
{\em Topology} {\bf 28} (1989), no. 2, 171--195.

\bibitem[PS1]{PS1}{\sc W.\ Pardon, M.\ Stern},
$L^2$-$\dq$-cohomology of complex projective varieties, {\em J. Amer. Math. Soc.} {\bf 4} (1991), no. 3, 603--621.

\bibitem[PS2]{PS2}{\sc W.\ Pardon, M.\ Stern},
Pure Hodge structure on the $L^2$-cohomology of varieties with isolated singularities,
{\em J. reine angew. Math.} {\bf 533} (2001), 55--80.

\bibitem[R1]{Rd}{\sc W.\ Rudin}, {\em Functional Analysis}, 
International Series in Pure and Applied Mathematics, McGraw-Hill, New York, 1991.

\bibitem[R2]{Rp1}{\sc J.\ Ruppenthal}, About the $\dq$-equation at isolated singularities with regular exceptional set,
{\em Internat. J. Math.} {\bf 20} (2009), no. 4, 459--489.

\bibitem[R3]{Rp7}{\sc J.\ Ruppenthal}, The $\dq$-equation on homogeneous varieties with an isolated singularity,
{\em Math. Z.} {\bf 263} (2009), 447--472.

\bibitem[R4]{Rp8}{\sc J.\ Ruppenthal},
$L^2$-theory for the $\dq$-operator on compact complex spaces, {\sf arXiv:1004.0396},
ESI-Preprint {\bf 2202}, submitted.

\bibitem[R5]{Rp9}{\sc J.\ Ruppenthal},
Compactness of the $\dq$-Neumann operator on singular complex spaces,
{\em J. Funct. Anal.} {\bf 260} (2011), no. 11, 3363-3403.


\bibitem[R6]{Rp10}{\sc J.\ Ruppenthal},
$L^2$-theory for the $\dq$-operator on complex spaces with isolated singularities, Preprint 2011, {\sf arXiv:1110:2373},
submitted.


\bibitem[RZ1]{RZ1} {\sc J.\ Ruppenthal, E.\ S.\ Zeron}, An explicit $\dq$-integration formula
for weighted homogeneous varieties, {\em Mich. Math. J.} {\bf 58} (2009), 321--337.

\bibitem[RZ2]{RZ2} {\sc J.\ Ruppenthal, E.\ S.\ Zeron}, An explicit $\dq$-integration formula
for weighted homogeneous varieties II, forms of higher degree, {\em Mich. Math. J.} {\bf 59} (2010), 283--295.

\bibitem[S1]{Sh}{\sc M.-C.\ Shaw}, Global solvability and regularity for $\dq$ on an annulus between two weakly pseudoconvex domains,
{\em Trans. Amer. Math. Soc.} {\bf 291} (1985), 255--267.


\bibitem[S2]{Siu00}{\sc Y.-T.\ Siu},
Analytic sheaf cohomology groups of dimension $n$ of $n$-dimensional non-compact complex manifolds,
{\em Pacific J. Math.} {\bf 28} (1969), 407--411.


\bibitem[S3]{Siu0}{\sc Y.-T.\ Siu},
Analyticity of sets associated to Lelong numbers and the extension of closed positive currents,
{\em Invent. Math.} {\bf 27} (1974), 53--156.

\bibitem[S4]{Siu}{\sc Y.-T.\ Siu},
Invariance of Plurigenera, {\em Invent. Math.} {\bf 134} (1998), no. 3, 661--673.

\bibitem[S5]{S}{\sc E.\ Straube},
{\em Lectures on the $L^2$-Sobolev theory of the $\overline{\partial}$-Neumann problem},
ESI Lectures in Mathematics and Physics, European Mathematical Society (EMS), Z\"urich, 2010.
\end{thebibliography}
\end{document}